\documentclass[titlepage,12pt]{amsart}
\usepackage{amssymb}

\newtheorem{Lemma}      {Lemma} [section]
\newtheorem{Theorem}    [Lemma] {Theorem}
\newtheorem{Definition} [Lemma] {Definition}
\newtheorem{Corollary}  [Lemma] {Corollary}
\newtheorem{Proposition}[Lemma] {Proposition}

\newcommand{\ov} [1]{{\overline{#1}}}
\newcommand{\op}    {\operatorname}
\newcommand{\invol} {{\mathfrak I}}
\newcommand{\cclass}{{\mathcal Cl}}

\begin{document}
\title[Real Subpairs and FS-Indicators]{Real Subpairs and Frobenius-Schur Indicators of Characters in $2$-Blocks}
\author{John Murray}
\date{November 7, 2008.}
\address{Mathematics Department, National University of Ireland,\newline
	 Maynooth, Co. Kildare, Ireland.}
\email{John.Murray@maths.nuim.ie}

\begin{abstract}
Let $B$ be a real $2$-block of a finite group $G$. A defect couple of $B$ is a certain pair $(D,E)$ of $2$-subgroups of $G$, such that $D$ a defect group of $B$, and $D\leq E$. The block $B$ is principal if $E=D$; otherwise $[E:D]=2$. We show that $(D,E)$ determines which $B$-subpairs are real.

The involution module of\/ $G$ arises from the conjugation action of\/ $G$ on its involutions. We outline how $(D,E)$ influences the vertices of components of the involution module that belong to $B$.

These results allow us to enumerate the Frobenius-Schur indicators of the irreducible characters in $B$, when $B$ has a dihedral defect group. The answer depends both on the decomposition matrix of\/ $B$ and on a defect couple of\/ $B$. We also determine the vertices of the components of the involution module of\/ $B$.
\end{abstract}
\maketitle
\thispagestyle{empty}

\section{Introduction}\label{S:Introduction}

Throughout this paper $G$ is a finite group. We adopt the standard notation and results for the representation theory of\/ $G$, as expounded in \cite{NagaoTsushima}. In particular ${\mathcal O}$ is a complete discrete valuation ring  of characteristic $0$ with residue field $k:={\mathcal O}/J({\mathcal O})$ of characteristic $p>0$. Mostly, but not exclusively, $p=2$. We assume that $\op{frak}(\mathcal O)$ and $k$ are splitting fields for all subgroups of\/ $G$. We write ${\op{N}(X)\!=\!\op{N}_G(X)}$ for the normalizer and ${\op{C}(X)\!=\!\op{C}_G(X)}$ for the centralizer of a subset $X$ of a $G$-set. If\/ $n$ is a natural number, we use $\nu(n)$ for the highest power of\/ $2$ that divides $n$, and ${\mathfrak S}_n$ to denote the symmetric group of degree $n$.

Let $B$ be a $p$-block of\/ $G$. There is a corresponding primitive idempotent $e_B$ of the centre $\mathord{Z}(kG)$ of the group algebra $kG$ and a $k$-algebra map $\omega_B:\mathord{Z}(kG)\rightarrow k$ such that $\omega_B(e_B)=1_k$. Attached to $B$ is a set $\op{Irr}(B)$ of irreducible characters and a set $\op{IBr}(B)$ of irreducible Brauer characters of\/ $G$. As is standard, $k(B)=|\op{Irr}(B)|$ and $l(B)=|\op{IBr}(B)|$.

The notion of vertex of a module, admissible block and induced block are as in \cite{NagaoTsushima}. Suppose that $M$ is an indecomposable module for\/ $H\leq G$ and that $\op{C}(V)\leq H$, for some vertex $V$ of\/ $M$. Then $M$ belongs to an admissible $p$-block $b$ of\/ $H$. {\em Green's vertex theorem} \cite[(3.7a)]{Green} states that if\/ $L$ is an indecomposable $kG$-module such that $M$ is a component of the restriction $L{\downarrow_H}$ then $L$ belongs to $b^G$.

A {\em $B$-subsection} is a pair $(x,b)$, where $x$ is a $p$-element of\/ $G$ and $b$ is a block of the centralizer $\op{C}(x)$ of\/ $x$ in $G$ such that $b^G=B$. If\/ $\theta\in\op{IBr}(b)$, we call $(x,\theta)$ a {\em column} of\/ $B$. For $\chi\in\op{Irr}(B)$ and $y$ a $p'$-element in $\op{C}(x)$ we have $\chi(xy)=\sum_{(x,b)}\sum_{\theta\in\op{IBr}(b)}d_{\chi\theta}^{(x)}\theta(y)$, where $d_{\chi\theta}^{(x)}$ are algebraic integers called {\em generalized decomposition numbers}. The group $G$ acts by conjugation on subsections and columns; $(x,\theta)^g=(x^g,\theta^g)$, for all $g\in G$. We say that $(x,\theta)$ is {\em real} if\/ $(x,\theta)^g=(x^{-1},\overline\theta)$, for some $g\in G$. Now identify conjugate $B$-subsections and columns. Brauer \cite{BrauerZurII} showed:
$$
k(B)=\sum_{(x,b)} l(b),\quad\mbox{where $(x,b)$ ranges over the $B$-subsections.}
$$

\begin{Lemma}\label{L:number_real_characters}
The number of real irreducible characters in $B$ equals the number of real columns in $B$.
\begin{proof}
Let $\chi\in\op{Irr}(B)$ and let $(x,\theta)$ be a column of\/ $B$.  Then
$$
d_{\chi\theta}^{(x)}=\sum\langle\chi\downarrow_{\op{C}_G(x)},\psi\rangle\frac{\psi(x)}{\psi(1)}d_{\psi\theta},
$$
where $\psi$ runs over the irreducible characters of\/ $\op{C}_G(x)$. Thus
$$
d_{\ov{\chi}\ov{\theta}}^{(x^{-1})}=\sum\langle\ov{\chi}\downarrow_{\op{C}_G(x)},\ov{\psi}\rangle\frac{\ov{\psi}(x^{-1})}{\ov{\psi}(1)}d_{\ov{\psi}\ov{\theta}}=d_{\chi\theta}^{(x)}.
$$

The generalized decomposition matrix $[d_{\chi\theta}^{(x)}]$ is nonsingular. It then follows from the previous paragraph and Brauer's permutation lemma that the number of\/ $\chi\in\op{Irr}(B)$ with $\chi=\ov{\chi}$ equals the number of Brauer subsections $(x,\theta)$ such that $(x,\theta)$ is conjugate to $(x^{-1},\ov\theta)$.
\end{proof}
\end{Lemma}

The {\em extended centralizer} of\/ $g\in G$ is ${\op{C}^*(g)}\!=\!{\op{C}_G^*(g)\!:=\!\op{N}(\{g,g^{-1}\})}$. So ${[\op{C}^*(g):\op{C}(g)]\leq2}$. We use ${\op{cl}(g)\!=\!\op{cl}_G(g)}$ to denote the $G$-conjugacy class of\/ $g$.

Now let $\op{char}(k)\!=\!2$ and let $B$ be a real $2$-block of\/ $G$. As noted by R. Gow \cite{Gow_JAlg}, there is a real conjugacy class ${\mathcal C}$ of\/ $G$, necessarily $2$-regular, so that ${\mathcal C}$ is contained in the support of\/ $e_B$, and $\omega_B({\mathcal C}^+)\ne0$. Here ${\mathcal C}^+=\sum_{c\in C}c$ in $\op{Z}(kG)$. Any such ${\mathcal C}$ is called a {\em real defect class} of\/ $B$.

\begin{Definition}\label{D:defect_couple}
A {\em defect couple} of a real $2$-block $B$ is a pair $(D,E)$, where $E$ is a Sylow $2$-subgroup of\/ $\op{C}^*(c)$, $D=\op{C}_E(c)$ is a Sylow $2$-subgroup of\/ $\op{C}(c)$, and $c$ is an element of a real defect class of\/ $B$.
\end{Definition}

As $D$ is a defect group of\/ $B$, it is unique up to $G$-conjugacy. Gow \cite{Gow_Osaka} showed that $E$ is also unique up to $G$-conjugacy. He referred to $E$ as an {\em extended defect group} of\/ $B$. In fact, the pair $(D,E)$ is unique up to $G$-conjugacy. This can be deduced from Proposition 14 in \cite{Murray_Osaka}.

Note that if\/ $\beta$ is an admissible real $2$-block of\/ $H\leq G$ with defect couple $(D_1,E_1)$, then by \cite[2.1]{Murray_GroupT} there is a defect group $D_2$ of\/ $\beta^G$ such that $D_1\leq D_2$ and $(D_2,D_2E_1)$ is a defect couple of\/ $\beta^G$.

The Frobenius-Schur (FS)-indicator $\epsilon(\chi)$ of a character $\chi$ of\/ $G$ is the average value of\/ $\chi(g^2)$, as $g$ ranges over $G$. If\/ $\chi$ is irreducible, it is known that $\epsilon(\chi)=0,+1,-1$, depending on whether $\chi$ is non-real, afforded by a real representation, or real-valued but not afforded by a real representation, respectively. The following is Corollary 2.5 in \cite{Murray_GroupT}:

\begin{Lemma}\label{L:2ElementSum}
Let $(D,E)$ be a defect couple of\/$B$ and let $x\in D$. Then
$$
\sum_{\chi\in\op{Irr}(B)}\epsilon(\chi)\chi(x)\geq 0,
$$
with $>0$ if and only if\/ ${E\!=\!D\langle e\rangle}$, where $e^2$ is $G$-conjugate to $x$.
\end{Lemma}

The columns of\/ $B$, weighted by FS-indicators, are locally determined, as Brauer has shown:

\begin{Lemma}\label{L:BrauerSum}
Let $(x,\theta)$ be a column of\/ $B$, with $\theta\in\op{Ibr}(b)$. Then
$$
\sum_{\chi\in\op{Irr}(B)}\epsilon(\chi)d_{\chi,\theta}^{(x)}\,=
\sum_{\psi\in\op{Irr}(b)}\epsilon(\psi)d_{\psi,\theta}^{(x)}.
$$
The value of these sums is non-negative if\/ $l(b)=1$.
\begin{proof}
The equality is Theorem (4A) of \cite{BrauerIII}. This holds even if\/ $B$ is a $p$-block, with $p$ odd. Suppose that $l(b)=1$. Then 
$$
\sum_{\psi\in\op{Irr}(b)}\epsilon(\psi)d_{\psi,\theta}^{(x)}\theta(1)=\sum_{\psi\in\op{Irr}(b)}\epsilon(\psi)\psi(x)\geq0,\quad\mbox{by Lemma \ref{L:2ElementSum}.}
$$
\end{proof}
\end{Lemma}

The following elementary lemma is useful for computing these sums:

\begin{Lemma}\label{L:SignedSum}
Let $m$ be an odd integer. Then for each integer $d>0$ with $|m|<2^d$, there is a {\em unique} $d$-tuple $(\varepsilon_0,\ldots,\varepsilon_{d-1})$ of signs such that
$$
\sum_{j=0}^{d-1}\varepsilon_j2^j\,=\,m.
$$
\end{Lemma}

We will use $[M:I]$ to denote the multiplicity with which an irreducible module $I$ occurs in a composition series of a module $M$. Then $[M:I]$ is well-defined, by the Krull-Schmidt theorem.

\begin{Lemma}\label{L:Indecomposable}
Let $M$ be a self-dual non-semisimple module and let $I$ be a self-dual irreducible submodule of\/ $M$. Then $[M:I]\geq2$.
\begin{proof}
We have $I\subseteq\op{soc}(M)\subseteq\op{rad}(M)$ and, by duality $I\subseteq M/\op{rad}(M)$.
\end{proof}
\end{Lemma}

The set of involutions in $G$, and the identity element, form the $G$-set
$$\invol=\invol_G:=\{g\in G\mid g^2=1_G\},$$
under conjugation. Suppose that $\op{char}(k)=2$. Then we call the permutation module $k\invol$ the {\em involution module} of\/ $G$. The involution module of\/ $B$ is the sum $k\invol B$ of all submodules of\/ $k\invol$ that belong to $B$. We showed in \cite{Murray_Osaka} that $k\invol B\ne0$ if and only if\/ $B$ is {\em strongly real}\/ i.e. $B$ is real and has a defect couple $(D,D\langle t\rangle)$, for some $t\in\invol$.

If\/ ${N\unlhd G}$ and ${X\!\subseteq\! G}$, we use ${\ov{X} =\!\{\!Nx\!\mid\! x\in X\}}$ to denote the image of\/ $X$ in the quotient group $\ov{G}:=G/N$.

\begin{Lemma}\label{L:Dominates}
Let $Q\unlhd G$ be a $2$-group and let $(D,E)$ be a defect couple of\/ $B$. Then there is a real $2$-block $\ov{B}$ of\/ $\ov{G}$, that is dominated by $B$, and that has defect couple $(\ov{D},\ov{E})$.
\begin{proof}
Write $\ov{x}\in k\ov{G}$ for the image of\/ $x\in kG$, under the natural $k$-algebra projection $kG\rightarrow k\ov{G}$. Let $\ov{B}_1,\ldots,\ov{B}_r$ be the blocks of\/ $\ov{G}$ that are dominated by $B$. The dual blocks $\ov{B}_i^o$ are also dominated by $B$.

Let ${\mathcal C}$ be a real defect class of\/ $B$. Then $\ov{{\mathcal C}}$ is a conjugacy class of\/ $\ov{G}$ and $\ov{{\mathcal C}}^{\,+}$ appears in the support of\/ $\ov{e_B}=\sum e_{\ov{B}_i}$, but not in the support of\/ $e_{\ov{B}_i}+e_{\ov{B}_{i}}^o$, for any $i$. It follows that we may choose $i$ such that $\ov{B}:=\ov{B}_i$ is real and $\ov{{\mathcal C}}^+$ appears in the support of\/ $e_{\ov{B}}$. As moreover $\omega_{\ov{B}}(\ov{{\mathcal C}}^+)=\omega_B({\mathcal C}^+)\ne 0$, we deduce that $\ov{{\mathcal C}}$ is a real defect class of\/ $\ov{B}$. Choose $c\in {\mathcal C}$ such that $E$ is a Sylow $2$-subgroup of\/ $\op{C}^*(c)$ and $D=E\cap\op{C}(c)$. Then $c$ is the unique element of odd order in $\ov{c}$. So $\ov{\op{C}_G^*(c)}=\op{C}^*_{\ov{G}}(\ov{c})$. It follows that $(\ov{D},\ov{E})$ is a defect couple of\/ $\ov{B}$.
\end{proof}
\end{Lemma}

In the context of the lemma, it is known that if\/ $G=Q\op{C}(Q)$, then $\ov{B}$ is the only block of\/ $\ov{G}$ dominated by $B$.

Suppose that $N\unlhd G$ and $b$ is a block of\/ $N$. Then $\op{N}(b)$ denotes the stabilizer of\/ $b$ in $G$, and $\op{N}^*(b)$ the stabilizer of\/ $\{b,b^o\}$ in $G$, all under the conjugation action of\/ $G$ on the blocks of\/ $N$. So $\op{N}(b)\leq\op{N}^*(b)$ and $[\op{N}^*(b):\op{N}(b)]\leq 2$. We call $\op{N}^*(b)$ the {\em extended stabilizer} of\/ $b$.

\begin{Lemma}\label{L:real_normal}
Let $N\unlhd G$ and let $b$ be an admissible $2$-block of\/ $N$ such that $b^G$ is real and $G=\op{N}^*(b)$. Let $(D,E)$ be a defect couple of\/ $b^G$.
\begin{itemize}
\item[(i)] If\/ $b=b^o$ then $(D\cap N,E\cap N)$ is a defect couple of\/ $b$.
\item[(ii)] If\/ $b\ne b^o$, then $G=E\op{N}(b)$ and $D\cap N$ is a defect group of\/ $b$.
\item[(iii)] Let $F\leq E$ with $E={DF}$. Then $(D\cap{NF},E\cap{NF})$ is contained in a defect couple of the real block $b^{NF}$.
\end{itemize}
\begin{proof}
Write $e_b=\sum_{n\in N}\beta_nn$, where $\beta_n\in k$, for each $n$. As $b$ is admissible, each real defect class ${\mathcal C}$ of\/ $b^G$ is contained in $N$. Let $c\in{\mathcal C}$ have defect couple $(D,E)$ in $G$.

Assume the hypothesis of (i). Then $b$ is $G$-invariant and hence $e_b$ is the primitive idempotent in $\mathord{Z}(kG)$ corresponding to $b^G$. Decompose ${\mathcal C}={\mathcal C}_1\cup\ldots\cup {\mathcal C}_t$ into a union of conjugacy classes of\/ $N$, where $c\in {\mathcal C}_1$. Now $G$ permutes the ${\mathcal C}_i$ transitively among themselves. It follows that $\omega_b({\mathcal C}_i^+)=\omega_b({\mathcal C}_1^+)$, for all $i$. So $0\ne\omega_{b^G}({\mathcal C}^+)=t\omega_b({\mathcal C}_1^+)$. In particular $\omega_b({\mathcal C}_1^+)\ne0$. So ${\mathcal C}_1$ is a real defect class of\/ $b$. Then (i) follows from the fact that $E\cap N$ is a Sylow $2$-subgroup of\/ $\op{C}_N(c)$ and $D\cap N=(E\cap N)\cap\op{C}_N(c)$.

Suppose in addition the hypothesis of (iii). Now $F\leq\op{C}^*(c)$ and ${\mathcal C}_1$ is a real class of\/ $N$. It follows that ${\mathcal C}_1$ is a real conjugacy class of\/ $NF$. But $e_b$ is the block idempotent of\/ $b^{NF}$. So ${\mathcal C}_1$ is a real defect class of\/ $b^{NF}$. As $(D\cap NF, E\cap NF)$ is contained in a defect couple of\/ $c$ in $NF$, it is also contained in a defect couple of\/ $b^{NF}$.

Assume the hypothesis of (ii). Then $[G:\op{N}(b)]=2$ and $e_b+e_b^o$ is the primitive idempotent in $\mathord{Z}(kG)$ corresponding to $b^G$. Regard $e_b$ as the block idempotent of a non-real block of\/ $\op{N}(b)$. As $c$ appears with nonzero multiplicity in $e_b+e_b^o$, we have $\beta_c+\beta_{c^{-1}}\ne 0$. But $\beta$ is a class function of\/ $\op{N}(b)$. So $c$ is not conjugate to $c^{-1}$ in $\op{N}(b)$, and in particular $G=E\op{N}(b)$. It is a standard fact that $D\cap N$ is a defect group of\/ $b$.

Suppose in addition the hypothesis of (iii). Then $b$ and $b^o$ are conjugate in $NF$. So $b^{NF}$ is a real block. We may write ${\mathcal C}=({\mathcal C}_1\cup\ldots\cup {\mathcal C}_t)\cup( {\mathcal C}_1^o\cup\ldots\cup {\mathcal C}_t^o)$ into a union of non-real conjugacy classes of\/ $N$, where ${\mathcal C}_1\cup\ldots\cup {\mathcal C}_t$ is a single conjugacy class of\/ $\op{N}(b)$. As before $\omega_b({\mathcal C}_i^+)=\omega_b({\mathcal C}_1^+)$ and $\omega_b({\mathcal C}_i^{o+})=\omega_b({\mathcal C}_1^{o+})$, for all $i$. Thus $0\ne\omega_{b^G}({\mathcal C}^+)=t\omega_b(({\mathcal C}_1\cup {\mathcal C}_1^{o})^+)$. Now ${\mathcal C}_1\cup {\mathcal C}_1^o$ is a real conjugacy class of\/ $NF$, and $e_b+e_b^o$ is the block idempotent of\/ $b^{NF}$. We deduce that ${\mathcal C}_1\!\cup {\mathcal C}_1^o$ is a real defect class of\/ $b^{NF}$\!. As $(D\cap NF, E\cap NF)$ is contained in a defect couple of\/ $c$ in $N\!F$\!, it is contained in a defect couple of\/ $b^{NF}$.
\end{proof}
\end{Lemma}

We now describe the rest of the paper. In Sections \ref{S:real_subpairs} and \ref{S:InvolutionModule} our only assumptions are that $p=2$, $B$ is a real $2$-block of\/ $G$, and $(D,E)$ is a defect couple of\/ $B$. We prove a number of results on the Brauer category of\/ $B$ in Section \ref{S:real_subpairs}. Our main result, Theorem \ref{T:real_subpairs}, shows how $(D,E)$ determines the reality of each $B$-subpair.

Suppose that $T$ is a conjugacy class of involutions in $G$. Let $kTB$ denote the sum of the components of the associated permutation module that belong to $B$. In Section \ref{S:InvolutionModule} we show how $(D,E)$ influences the vertices of the components of\/ $kTB$: Theorem \ref{T:C_D(t)} gives a connection between these vertices and the centralizers in $D$ of elements of\/ $T\cap E$.

From Section \ref{S:Extensions} onwards, we assume that $B$ is a real $2$-block that has a dihedral defect group $D$. We first enumerate the possible isomorphism types of an extension $E$ of\/ $D$ with $[E:D]=2$. Table \ref{Tb:E} describes the $E$-conjugacy classes in $E\backslash D$ for each of these types.

In Section \ref{S:local}, we determine, for each type of\/ $E$, the real subpairs of\/ $B$ (Lemma \ref{L:real_subpairs}), and the real irreducible characters in $B$ that have positive height (Theorem \ref{T:real_characters}).

Section \ref{S:Morita} reviews known results on the Morita equivalence classes of\/ $2$-blocks with dihedral defect groups. In particular, Lemma \ref{L:MoritaTypes} partially describes the generalized decomposition matrices.

In Sections \ref{S:TotallySplit} through \ref{S:Type(e)} we determine the FS-indicators of the irreducible characters in $B$, and the vertices of the components of\/ $k\invol B$, according to the isomorphism type of the extended defect group $E$, and the Morita equivalence class of\/ $B$.

Our results for blocks with dihedral defect group are summarized in Table \ref{Tb:main}. We believe that the methods of this paper could be applied to all blocks with tame representation type (i.e. dihedral, semidihedral or quaternion defect group) in order to determine the Frobenius-Schur indicators of the irreducible characters and the vertices of components of the involution module.

\section{Real subpairs}\label{S:real_subpairs}

The purpose of this section is to show how a defect couple of a real $2$-block determines the reality of its subpairs. We assume that $B$ is a real $2$-block of\/ $G$, and let $(D,E)$ be a defect couple of\/ $B$.

A {\em $B$-subpair} consists of\/ $(Q,b_Q)$ where $Q$ is a $2$-subgroup of\/ $G$ and $b_Q$ is a $2$-block of\/ $Q\op{C}(Q)$ such that $b_Q^G=B$. By abuse of notation $\op{N}(b_Q)$ and $\op{N}^*(b_Q)$ denote the stabilizer and extended stabilizer of\/ $b_Q$ in $\op{N}(Q)$. A $B$-subpair of the form $(D,b_D)$ is known as a {\em Sylow $B$-subpair} (or a {\em root} of\/ $B$). Brauer showed that all Sylow $B$-subpairs $(D,b_D)$ are $\op{N}(D)$-conjugate, and also that $[\op{N}(b_D):\op{C}(D)D]$ is odd.

The set of\/ $B$-subpairs is a poset. Containment $\leq$ is generated by normal containment: $(Q,b_Q)\unlhd(R,b_R)$ if\/ $R\leq\op{N}(b_Q)$ and $b_R^{R\op{C}(Q)}=b_Q^{R\op{C}(Q)}$. Given a $B$-subpair $(Q,b_Q)$, there is a Sylow $B$-subpair $(D,b_D)$ such that $(Q,b_Q)\leq(D,b_D)$. For each Sylow $B$-subpair $(D,b_D)$ and $Q\leq D$, there is a unique $B$-subpair $(Q,b_Q)$ such that $(Q,b_Q)\leq(D,b_D)$. We will refer to this as the {\em uniqueness} property of subpairs.

\begin{Lemma}\label{L:going_down}
Suppose that $(Q,b_Q)$ and $(R,b_R)$ are $B$-subpairs such that $b_R$ is real and $(Q,b_Q)\leq(R,b_R)$. Then $b_Q$ is real.
\begin{proof}
We have $(Q,b_Q)\leq(R,b_R)$, and $(Q,b_Q^o)\leq(R,b_R^o)=(R,b_R)$. So $b_Q^o=b_Q$, by the uniqueness property of subpairs.
\end{proof}
\end{Lemma}

We now examine the relationship between the extended defect groups of\/ $B$ and the extended stabilizers of its  Brauer subpairs.

\begin{Lemma}\label{L:inertia*} The following hold:
\begin{itemize}
\item[(i)] If\/ $(D,b_D)$ is a Sylow $B$-subpair, then there is a defect couple $(D,E)$ of\/ $B$ such that
$b_D^{E\op{C}(D)}$ is real with defect couple $(D,E)$.
\item[(ii)] If $(D,E)$ is a defect couple of $B$ then there is a Sylow $B$-subpair $(D,b_D)$ such that $b_D^{E\op{C}(D)}$ is real with defect couple $(D,E)$.
\end{itemize}
\begin{proof}
Suppose first that $b_D$ is real, with defect couple $(D,E)$. Then $(D,E)$ is a defect couple of\/ $B$. Moreover, $E\op{C}(D)=D\op{C}(D)$. Suppose then that $b_D$ is not real. As $(D,b_D)$ and $(D,b_D^o)$ are Sylow $B$-subpair, they are conjugate in $\op{N}(D)$, and hence in $\op{N}^*(b_D)$. Set $b=b_D^{\op{N}^*(b_D)}$. Then each defect couple $(D,E)$ of\/ $b$ is a defect couple of\/ $B$. Lemma \ref{L:real_normal} implies that $b_D^{E\op{C}(D)}$ is real with defect couple $(D,E)$.

Now let $(D,E)$ be a defect couple for $B$, and let $(D,b_D)$ be any Sylow $B$-subpair. By part (i) there exists a defect couple $(D,E_1)$ for $B$ such that $b_D^{E\op{C}(D)}$ is real with defect couple $(D,E_1)$. But the defect couples of\/ $B$ are conjugate in $G$. So there exists $x\in G$ such that $(D,E)=(D,E_1)^x$. Now replace $b_D^x$ by $b_D$. Then $(D,b_D)$ a Sylow $B$-subpair such that $b_D^{E\op{C}(D)}$ is real with defect couple $(D,E)$.
\end{proof}
\end{Lemma}

A {\em $B$-subtriple} consists of $(Q,b_Q,R)$ where:
\begin{itemize}
\item $(Q,b_Q)$ is a $B$-subpair;
\item $R\leq\op{N}^*(b_Q)$ and $b_Q^*:=b_Q^{R\op{C}(Q)}$ is real;
\item $R$ is an extended defect group of\/ $b_Q^*$.
\end{itemize}
We say that this is a {\em Sylow $B$-subtriple} if\/ $Q$ is a defect group of\/ $B$. For two $B$-subtriples an inclusion $(Q,b_Q,R)\leq(S,b_S,T)$ occurs if\/ $(Q,b_Q)\leq(S,b_S)$, in the sense of\/ $B$-subpairs, and  $T=RS$. We write $(Q,b_Q,R)\unlhd(S,b_S,T)$ if in addition $Q\unlhd S$. Note that if\/ $(Q,b_Q)$ is a real $B$-subpair and if\/ $R$ is an extended defect group of\/ $b_Q$, then $b_Q^*=b_Q$ and $(Q,b_Q,R)$ is a $B$-subtriple.

\begin{Lemma}\label{L:eo_action}
Suppose that $(D,b_D,E)$ is a Sylow $B$-subtriple and that ${Q\leq D}$. Then $b_Q^o=b_{Q^e}^{e^{-1}}$, for all $e\in E$ such that $E=D\langle e\rangle$.
\begin{proof}
We have $D^e=D$ and $b_D^e=b_D^o$. So $(Q^e,b_Q^{eo})\leq(D^e,b_D^{eo})=(D,b_D)$. Then $b_{Q^e}=b_Q^{eo}$, by the uniqueness property of subpairs. The lemma follows from this.
\end{proof}
\end{Lemma}

\begin{Lemma}\label{L:Sylow_subtriple_conjugacy}
All Sylow $B$-subtriples are conjugate in $G$.
\begin{proof}
Let $(D_1,b_1,E_1)$ and $(D_2,b_2,E_2)$ be Sylow $B$-subtriples. Then $(D_i,b_i)$ are Sylow $B$-subpairs. So $(D_2,b_2)^x=(D_1,b_1)$, for some $x\in G$. Now $(D_1,b_1,E_2^x)$ is a Sylow $B$-subtriple. In particular $(D_1,E_2^x)$ is a defect couple of\/ $b_1^*$ regarded as a block of\/ $\op{N}^*(b_1)$. But $(D_1,E_1)$ is another defect couple of this block. So there exists $y\in\op{N}^*(b_1)$ such that $(D_1,E_2^x)^y=(D_1,E_1)$. Now $b_1^y=b_1$ or $b_1^o$. In the former case set $z:=1$. Suppose we have the latter case. Then there exists $z\in E_1\op{C}(D_1)$ such that $(b_1^o)^z=b_1$. Then $(D_2,b_2,E_2)^g=(D_1,b_1,E_1)$, with $g=xyz$.
\end{proof}
\end{Lemma}

Let $(Q,b_Q,R)$ be a $B$-subtriple. So $b_Q^*$ is a real block of\/ $R\op{C}(Q)$. We can choose a defect couple $(S,T)$ of\/ $b_Q^*$, regarded as a block of\/ $\op{N}^*(b_Q)$, such that $T=RS$. Then $S$ is a defect group of\/ $b_Q$, regarded as a block of\/ $\op{N}(b_Q)$. Lemma \ref{L:inertia*} implies that there is a Sylow $b_Q^*$-subpair $(S,b_S)$ in $N^*(b_Q)$ such that $b_S^{T\op{C}(S)}$ is a real block with defect couple $(S,T)$. Moreover, $b_S^{S\op{C}(Q)}$ covers $b_Q$ or $b_Q^o$, so by switching to $b_S^o$, if necessary, we may assume that $(Q,b_Q)\leq(S,b_S)$. It follows that $(S,b_S,T)$ is a $B$-subtriple and $(Q,b_Q,R)\leq(S,b_S,T)$. We call any such $B$-subtriple a {\em normalizer-subtriple} of\/ $(Q,b_Q,R)$.

Set $(Q_1,b_1,R_1)\!:=\!(Q,b_Q,R)$ and inductively choose $(Q_{i+1},b_{i+1},R_{i+1})$ to be a normalizer subtriple of\/ $(Q_i,b_i,R_i)$, for $i=1,2,\ldots$. Since $(Q_{i+1},b_{i+1})$ is a normaizer-subpair of\/ $(Q_i,b_i)$, it follows that the sequence $(Q_1,b_1,R_1)\unlhd(Q_2,b_2,R_2)\unlhd(Q_3,b_3,R_3)\unlhd\ldots$ terminates at a Sylow $B$-subtriple $(D,b_D,E)$, where $(D,E)$ is a defect couple of\/ $B$. In this case $E=R_1D$. We say that $Q$ is {\em really well-placed in $D$} with respect to $(D,b_D,E)$ if such a sequence exists.

\begin{Lemma}\label{L:well-placed}
Let $(D,b_D,E)$ be a Sylow $B$-subtriple and let $(Q,b_Q,R)$ be any $B$-subtriple such that $(Q,b_Q)\leq(D,b_D)$. Then there exists $g\in G$ such that $Q^g$ is really well-placed in $D$ with respect to $(D,b_D,E)$ and $(Q,b_Q,R)^g\leq(D,b_D,E)$.
\begin{proof}
By the discussion above, there is a Sylow $B$-subtriple $(D_1,b_1,E_1)$ that contains $(Q,b_Q,R)$, such that $Q$ is {\em really well-placed in $D_1$} with respect to $(D_1,b_1,E_1)$. By Lemma \ref{L:Sylow_subtriple_conjugacy} there exists $g\in G$ such that ${(D_1,b_1,E_1)^g\!=\!(D,b_D,E)}$. The lemma follows from this.
\end{proof}
\end{Lemma}

Our next Lemma is concerned with strongly real $B$-subpairs.

\begin{Lemma}\label{L:strongly_real}
Let $(D,b_D,E)$ be a Sylow $B$-subtriple and let $(Q,b_Q)\leq(D,b_D)$ be such that $E=D\op{C}_E(Q)$. Then $b_Q$ is real and has a defect couple containing $(Q\op{C}_D(Q),Q\op{C}_E(Q))$. In particular, if\/ $E=D\langle t\rangle$, for some $t\in\invol_{\op{C}(Q)}$, then $b_Q$ is strongly real.
\begin{proof}
If\/ $Q=D$, the result follows from Lemma \ref{L:real_normal}. So we may assume that $Q<D$. Set $b_Q,D_0:=Q$ and inductively $D_i:=\op{N}_D(D_{i-1})$ and $E_i:=\op{N}_E(D_{i-1})$, for $i>0$. Then $[E_i:D_i]\leq 2$ and $E=DE_i$. For each positive $i$ there is a unique subpair $(Q,b_Q)\leq(D_i,b_i)\leq(D,b)$. Let $t$ be the smallest positive integer such that $D_t=D$ and $E_t=E$.

We use downwards induction on $i$ to prove that $b_i^*\!:=\!b_i^{E_i\op{C}(D_i)}$ is real with a defect couple containing $(D_i,E_i)$. The base case $i=t$ follows from the hypothesis. Assume the result for $i+1$. Set ${b_{i+1}^{**}\!:=b_{i+1}^{E_{i+1}{\rm C}(D_i)}}$. Then ${b_{i+1}^{**}={(b_{i+1}^*)^{E_{i+1}{\rm C}(D_i)}}}$. So $b_{i+1}^{**}$ is real with a defect couple containing $(D_{i+1},E_{i+1})$. Part (iii) of Lemma \ref{L:real_normal}, applied with $N=D_i\op{C}(D_i)$, $b=b_i$, $G=E_{i+1}\op{C}(D_i)$ and $F=E_i$, gives the inductive step.

The previous paragraph shows that $b_1^*$ is a real block of\/ $\op{N}_E(Q)\op{C}(Q)$ that covers the block $b_Q$ of\/ $Q\op{C}(Q)$. We apply part (iii) of Lemma \ref{L:real_normal}, with $N=Q\op{C}(Q)$, $b=b_Q$, $G=\op{N}_E(Q)\op{C}(Q)$ and $F=\op{C}_E(Q)$. Then $NF=Q\op{C}(Q)$. So $b_Q$ is a real block with a defect 
couple containing $(Q\op{C}_D(Q),Q\op{C}_E(Q))$.
\end{proof}
\end{Lemma}

We use $\sim$ for $G$-conjugacy. Our main result on real subpairs is:

\begin{Theorem}\label{T:real_subpairs}
Let $(D,b_D,E)$ be a Sylow $B$-subtriple. Set
$$
{\mathcal R}(D,b_D,E):=\{(Q,b_Q)\leq(D,b_D)\mid E=D\op{C}_E(Q)\}.
$$
Then a $B$-subpair $(R,b_R)$ is real if and only if\/ $(R,b_R)\sim(Q,b_Q)$, for some $(Q,b_Q)\in{\mathcal R}(D,b_D,E)$. Moreover, $(R,b_R)$ is strongly real if and only if\/ $(R,b_R)\sim(Q,b_Q)$, where $E=D\langle t\rangle$, for some $t\in\invol_{\op{C}(Q)}$.
\begin{proof}
Suppose first that $(Q,b_Q)\in{\mathcal R}(D,b_D,E)$. Then by Lemma \ref{L:inertia*} we can choose $e\in\op{C}_E(Q)$ such that $b_D^{oe}=b_D$. Then $(Q,b_Q^o)=(Q,b_Q^o)^e\leq(D,b_D^o)^e=(D,b_D)$. But $(Q,b_Q)$ is the unique $B$-subpair involving $Q$ and contained in $(D,b_D)$. So $(Q,b_Q^o)=(Q,b_Q)$, whence $b_Q$ is real.
Now suppose that $x\in G$ and $(Q,b_Q)^x\leq(D,b_D)$. Then $(b_Q^x)^o=(b_Q^o)^x=b_Q^x$. So $(Q,b_Q)^x$ is a real $B$-subpair contained in $(D,b_D)$. If in addition $E=D\langle t\rangle$, for some $t\in\invol_{\op{C}(Q)}$, then $b_Q$ is strongly real, by Lemma \ref{L:strongly_real}. This completes the `if' part of the theorem.

The `only if' statement follows from Lemma \ref{L:well-placed}.
\end{proof}
\end{Theorem}

\section{The involution module}\label{S:InvolutionModule}

In this section we survey known facts about the involution module $k\invol$ of\/ $G$, and prove a number of new results about the vertices of its components. As in the previous section, $B$ is a real $2$-block of\/ $G$ that has defect couple $(D,E)$. Fix a conjugacy class $T\subseteq\invol$ of\/ $G$.

Suppose that $I$ is an irreducible $B$-module, with Brauer character $\theta$. The projective character of\/ $G$ associated with $\theta$ is $\Phi=\sum_{\chi\in\op{Irr}(B)} d_{\chi\theta}\,\chi$. Here the $d_{\chi\theta}$ are non-negative integers called {\em decomposition numbers}. Lemma 1 of \cite{Robinson} implies that
\begin{equation}\label{E:Smult}
\epsilon(\Phi)=[k\invol:I].
\end{equation}

Now suppose that $I$ is also self-dual and projective. Then $D=\langle1_G\rangle$ and $\Phi$ is the unique irreducible character in $B$. So $[k\invol:I]=\epsilon(\Phi)=+1$. Conversely, by \cite{Murray_JAlg}, each projective component of\/ $k\invol$ is self-dual and irreducible and belongs to a real $2$-block with a trivial defect group.

Let $M$ be a component of\/ $kTB$ and let $V$ be a vertex of\/ $M$. As $M$ has a trivial vertex, the Green correspondent $P$ of\/ $M$ with respect to $(G,V,\op{N}(V))$ is a component of\/ $k\op{C}_T(V)$. In particular, $P$ is projective as $\op{N}(V)/V$-module. Clifford theory shows that $P{\downarrow_{V\op{C}(V)}}=m\sum R^n$, where $m>0$, $R$ is indecomposable with vertex $V$ and $n$ ranges over a set of representatives for the cosets of the stabilizer of\/ $R$ in $\op{N}(V)$. Moreover, $R$ is projective as ${V{\rm C}(V)}/V$-module, and is a component of\/ $k\op{C}_T(V){\downarrow_{V\op{C}(V)}}$. Let $b_V$ be the block of\/ $V{\rm C}(V)$ that contains $R$. Then $(V,b_V)$ is a $B$-subpair and $P$ belongs to $b_V^{\op{N}(V)}$. Given $V$, $R$ is uniquely determined up to $\op{N}(V)$-conjugacy. We call $(V,b_V)$ a {\em vertex $B$-subpair} and $R$ a {\em $b_V$-root} of\/ $M$. These notions are set out in a more general context in \cite{Kulshammer}.

\begin{Lemma}\label{L:detectVertex}
Suppose that $D$ is non-trivial and $(V,b_V)$ is a vertex $B$-subpair of a component of\/ $kTB$. Then there is an involution $v\in V$, and a $B$ subsection $(v,b)$, such that some component of\/ $k\op{C}_T(v)b$ has vertex $b$-subpair $(V,b_V)$.
\begin{proof}
The hypothesis on $D$, and the discussion prior to the lemma, implies that $V\ne\{1_G\}$. Choose an involution $v\in\mathord{Z}(V)$. Then ${V\!\op{C}(V)}\leq\op{C}(v)$. Let $M$ be a component of\/ $kTB$ that has vertex $B$-subpair $(V,b_V)$. Then some component of\/ $M{\downarrow_{\op{C}(v)}}$ has vertex subpair $(V,b_V)$. Green's vertex theorem implies that this module belongs to a block $b$ such that $(v,b)$ is a $B$-subsection.
\end{proof}
\end{Lemma}

Our main result in this section sharpens Theorem 1.5 of \cite{Murray_GroupT}:

\begin{Theorem}\label{T:C_D(t)}
The following hold:
\begin{itemize}
\item[(i)] If\/ $kTB\ne 0$ and $M$ is a component of\/ $kTB$, then there exists $t\in T$ such that $E=D\langle t\rangle$ and $\op{C}_D(t)$ contains a vertex of\/ $M$.
\item[(ii)] Suppose that ${E\!=\!D\langle t\rangle}$, with $t\in T$, but $\op{C}_D(t)\!\not<_G\!\op{C}_D(s)$ for any $s\!\in{T\!\cap\!Dt}$. Then some component of\/ $kTB$ has vertex $\op{C}_D(t)$.
\end{itemize}
\begin{proof}[Proof of (i)]
For both parts we use induction on $|D|$. The base case $|D|=1$, is dealt with by Theorem 19 of \cite{Murray_Osaka}. So we may assume that $|D|>1$. Let $V$ be a vertex of\/ $M$ and set $H:=V\op{C}(V)$. Then $V\ne\{1_G\}$. Let $(V,b)$ be a vertex $B$-subpair and let $R$ be a $b$-root of\/ $M$. So $R$ is a component of\/ $k\op{C}_T(V)b$. Choose $t\in\op{C}_T(V)$ such that $R$ is a component of\/ $k\op{cl}_H(t)$.

Set $\ov{H}:=H/Z$, where $Z\leq\op{Z}(V)$ has order $2$. By Lemma \ref{L:Dominates}, there is a unique block $\ov{b}$ of\/ $\ov{H}$ that is dominated by $b$. Now $\op{N}_H(\ov{t})$ is the preimage of\/ $\op{C}_{\ov{H}}(\ov{t})$ in $H$ and $[\op{N}_H(\ov{t}):\op{C}_H(t)]\leq 2$. So either $k_{\op{N}_H(\ov{t})}{\uparrow^H}=k\op{cl}_H(t)$ or there is a short exact sequence of\/ $kH$-modules
$$
0\rightarrow k_{\op{N}_H(\ov{t})}{\uparrow^H}\rightarrow k\op{cl}_H(t)\rightarrow k_{\op{N}_H(\ov{t})}{\uparrow^H}\rightarrow 0.
$$
In any event each composition factor of\/ $R$ is a composition factor of\/ $k\op{cl}_{\ov{H}}(\ov{t})\ov{b}$. As $\ov{b}$ has a smaller defect group than $B$, the inductive hypothesis implies that $\ov{b}$ has a defect couple of the form $(\ov{D_1},\ov{D_1}\langle\ov{t}\rangle)$. Lemma \ref{L:Dominates} now shows that $b$ has a defect couple $(D_1,D_1\langle t\rangle)$. As $b^G=B$, we may assume that $D$ is chosen so that $D_1\leq D$ and $(D,D\langle t\rangle)$ is a defect couple of\/ $B$. Thus also $V\leq\op{C}_D(t)$.
\end{proof}
\begin{proof}[Proof of (ii)]
Set ${V:=\op{C}_D(t)}$ and ${H:=V\op{C}(V)}$. Let $\op{Br}_V\colon\op{Z}(kG)\rightarrow\op{Z}(k\op{C}(V))$ be the Brauer homomorphism with respect to $V$. Then $\op{Br}_V(e_B)=\sum e_b$, where $b$ ranges over the $2$-blocks of\/ $H$ such that $b^G=B$. Note that $(b^o)^G=b^G$, as $B$ is real. Choose $c$ in a real defect class of\/ $B$ such that $D\langle t\rangle$ is a Sylow $2$-subgroup of\/ $\op{C}^*(c)$ and $c^t=c^{-1}$. In particular $c,t\in H$. Set ${\mathcal C}:=\op{cl}_H(c)$. Then ${\mathcal C}^+$ is a real conjugacy class of\/ $H$ that appears in $\op{Br}_V(e_B)$. As $e_b+e_b^o$ is supported on the non-real classes of\/ $H$, for each block $b$ of\/ $H$, it follows that there exists a real block $b$ of\/ $H$ such that ${\mathcal C}^+$ appears in $e_b$ and $b^G=B$. Theorem 3.3 of \cite{Murray_CAlg} implies that $b$ has a defect couple of the form $(D_1,D_1\langle t\rangle)$.

Lemma \ref{L:Dominates} implies that $b$ dominates a real block $\ov{b}$ of\/ $\ov{H}:=H/V$ such that $\ov{b}$ has defect couple $(\ov{D_1},\ov{D_1}\langle\ov{t}\rangle)$. As $|\ov{D_1}|<|D_1|\leq|D|$, the inductive hypothesis implies that $k\cclass_{\ov{H}}(\ov{t})\ov{b}\ne0$. The inflation of this module to $H$ is a quotient module of\/ $k\cclass_H(t)b$. So $k\cclass_H(t)b\ne 0$.

Let $M_1$ be a component of\/ $k\cclass_H(t)b$. Then $V$ is contained in each vertex of\/ $M_1$. Let $M$ be a component of\/ $kT$ such that $M_1$ is a component of\/ $M{\downarrow_H}$. Then $M$ belongs to $B$ and $V$ is contained in some vertex $V_1$ of\/ $M$. Now from part (i) of this theorem $V_1\leq_G\op{C}_D(s)$, for some $s\in T\cap Dt$. Since $\op{C}_D(t)\leq V_1\leq_G\op{C}_D(s)$, it follows from the hypothesis that $V=V_1$ is a vertex of\/ $M$.
\end{proof}
\end{Theorem}

\begin{Corollary}
Suppose that $t$ is an involution in $\op{O}_2(G)$. Then every component of\/ $k_{\op{C}_G(t)}{\uparrow^G}$ belongs to the principal $2$-block of\/ $G$.
\begin{proof}
Let $B$ be a $2$-block of\/ $G$ such that $k_{\op{C}_G(t)}{\uparrow^G}B\ne 0$. Then each defect group $D$ of\/ $B$ contains $\op{O}_2(G)$. In particular $t\in D$. But then $D=D\langle t\rangle$ is an extended defect group of\/ $B$. It follows from this that $B$ has a real defect class in $\invol$. Since a defect class is necessarily $2$-regular, the class $\{1_G\}$ is a defect class of\/ $B$. As $B$ is real, it follows from this that $B$ is the principal $2$-block of\/ $G$.
\end{proof}
\end{Corollary}

We take the opportunity to correct an error in the hypothesis of Lemma 2.11 of \cite{Murray_GroupT}. We also strengthen the conclusion. Fortunately, we only used that lemma in its correct form.

\begin{Lemma}\label{L:Alperin}
Suppose that $D\unlhd G$ and $E=D\times\langle e\rangle$. Then there is a self-dual irreducible $B$-module $I$, such that $I$ has vertex $D$ and
$$
\bigoplus_{i\geq 0}I\otimes\op{rad}^i(k\invol_{\mathord{Z}(D)})/\op{rad}^{i+1}(k\invol_{\mathord{Z}(D)})
$$
is the sum of all components of\/ $k\invol B$ that have vertex $D$.
\begin{proof}
By hypothesis $e^2=1$ and $e\in\op{C}(D)$. This is what we needed, and used, in the proof of Lemma 2.11 of \cite{Murray_GroupT}. So that proof gives all but the last statement.

We make a very general remark. Let $H$ be a finite group, and let $M$ be an indecomposable $kH$-module that has vertex $V$ and source $k_V$. Set $f{M}$ as the Green correspondent of\/ $M$ with respect to $(H,V,\op{N}_H(V))$. Suppose that $V\leq W\leq\op{N}_H(V)$. Then $f{M}{\downarrow_W}$ is the sum of all components of\/ $M{\downarrow_W}$ that have vertex $V$.

We can now prove the last statement of the lemma. Adopt the notation of Lemma 2.11 of \cite{Murray_GroupT}. We apply the previous paragraph with $H=G\wr\Sigma$, $M=B$, $V=\Delta D\times\Sigma$ and $W=\Delta G\times\Sigma$. Then by Lemma 2.9 of \cite{Murray_GroupT}, $B(V){\downarrow_{\Delta G}}$ can be identified with the sum of all components of\/ $k\invol B^{\op{Fr}}$ that have vertex $D$. We get the equivalent result for $k\invol B$ using Lemma 2.8 of \cite{Murray_GroupT}.
\end{proof}
\end{Lemma}

\section{Extensions of dihedral $2$-groups}\label{S:Extensions}

We need to describe the degree $2$-extensions of $D$, where $D$ is a dihedral group of order $2^d$. Fix a presentation $D=\langle s,t\mid s^{2^{d-1}},t^2,(st)^2\rangle$. We shall adopt the following notation for subgroups of\/ $D$. The maximal cyclic subgroup $S:=\langle s\rangle$ of\/ $D$ has order $2^{d-1}$ and the centre $\op{Z}(D)=\langle s_1\rangle$ has order $2$. Also ${X_1:=\langle t\rangle}$ and ${Y_1:=\langle st\rangle}$ are subgroups of order $2$. Set $s_i:=s^{2^{d-1-i}}$, for ${i=1,\dots,d-1}$. Then $S_i:=\langle s_i\rangle$ is a cyclic group of order $2^i$. Moreover, $X_2:=\langle s_1,t\rangle$ and $Y_2:=\langle s_1,st\rangle$ are Klein-four groups, while $X_i:=\langle s_{i-1},t\rangle$ and $Y_i:=\langle s_{i-1},st\rangle$ are dihedral groups of order $2^i$, for $i\geq3$. In particular, $X_{d-1}$ and $Y_{d-1}$ are maximal subgroups of\/ $G$ that are dihedral of order $2^{d-1}$.

\begin{Proposition}\label{P:Dihedral_extensions}
There are $4$ isomorphism classes of degree $2$ extensions of\/ $D_8$. There are $5$ isomorphism classes of degree $2$ extensions $E=D\langle e\rangle$ of\/ $D\cong D_{2^d}$, for $d\geq 4$. These are:
\begin{itemize}
\item[(a)] $E=D\times\langle e\rangle$ where $e\in\op{C}(D)$ and $e^2=1$;
\item[(b)] $E=D\langle e\rangle$ where $e\in\op{C}(D)$ and $e^2=s_1$;
\item[(c)] $E=D_{2^{d+1}}$, a dihedral group of order $2^{d+1}$;
\item[(d)] $E=SD_{2^{d+1}}$, a semi-dihedral group of order $2^{d+1}$;
\item[(e)] $E=D\rtimes\langle e\rangle$ where $e^2=1$, $s^e=s_1s,t^e=t$;
\end{itemize}
\begin{proof}
Set $q=2^{d-1}$. The additive group of the ring ${\mathbb Z}_q$ of integers modulo $q$ is cyclic, while its group of units ${\mathbb U}_q$ is isomorphic to the direct product ${\mathbb Z}_2\times{\mathbb Z}_{2^{d-3}}$. The ${\mathbb Z}_2$ factor is generated by $-1$ and the involutions in ${\mathbb U}_q$ are $1,-1,2^{d-2}-1,2^{d-2}+1$. Form the semidirect product ${\mathbb U}_q\ltimes{\mathbb Z}_q$. So ${\mathbb U}_q$ acts by multiplication on ${\mathbb Z}_q$. We define an action of\/ ${\mathbb U}_q\ltimes{\mathbb Z}_q$ on $D$ via
$$
s^{(a,b)}:=s^a,\quad t^{(a,b)}:=s^{-b}t,\quad\mbox{for $a\in{\mathbb U}_q,b\in{\mathbb Z}_q$.}
$$
In this way ${\mathbb U}_q\ltimes{\mathbb Z}_q$ can be identified with the automorphism group of\/ $D$. The group of inner automorphisms of\/ $D$ is $I:=\langle-1\rangle\ltimes 2{\mathbb Z}_q$. So the outer automorphism group is isomorphic to ${\mathbb Z}_{2^{d-3}}\times{\mathbb Z}_{2}$. It follows that the image of\/ $E$ in the automorphism group of\/ $D$ is generated, modulo $I$ by one of\/ $(1,0),(-1,1),(2^{d-2}+1,0)$ or $(2^{d-2}-1,1)$. We set the image of\/ $e$ to be one of these four elements, in turn. The action of\/ $e$ on $D$ determines the coset $e^2\mathord{Z}(D)$. As $\mathord{Z}(D)$ has order $2$, we get two possibilities for $E$ in each case.

Suppose that the image of\/ $e$ is $(1,0)$. Then $e^2\in\op{Z}(D)$ and $E$ is the group of type (a) if\/ $e^2=1$, or the group of type (b) if\/ $e^2=s_1$.

Suppose that the image of\/ $e$ is $(-1,1)$. Then $e^2\in\{1,s_1\}$ and $E$ is a dihedral group if\/ $e^2=1$, or a semidihedral group, if\/ $e^2=s_1$.

From now on assume that $d\geq 4$. Suppose that the image of\/ $e$ is $(2^{d-2}+1,0)$. Then $e^2=1$ or $s_1$. In either case we obtain a group which is isomorphic to one of type (e).

We claim that the image of\/ $e$ cannot equal $(2^{d-2}-1,1)$. For suppose otherwise. Then $s^e=s_1s^{-1}$. So $e$ inverts $s_2\in\langle s^2\rangle$. On the other hand $e^2=s_2$ or $s_2^{-1}$. So $e$ centralizes $s_2$. This contradiction proves our claim.
\end{proof}
\end{Proposition}

\begin{Corollary}
${E\!:\!D}$ is non-split if and only if\/ $E$ is semi-dihedral. In all cases $E/D'$ splits over $D/D'$.
\begin{proof}
The first statement follows from the third column of Table \ref{Tb:E}, below. For $E\backslash D$ contains no involution if and only if\/ $E$ is of type (d). The second statement follows from the fact that if\/ $E$ is semi-dihedral then $D/D'\cong{\mathbb Z}_2^2$ and $E/D'\cong D_8$.
\end{proof}
\end{Corollary}

\begin{Corollary}\label{C:heightzero_FS+1}
Let $B$ be a $2$-block with a dihedral defect group. Then all real height zero irreducible characters in $B$ have FS-indicator $+1$.
\begin{proof}
This follows from Theorem 5.6 of \cite{Gow_JAlg} and the Corollary above.
\end{proof}
\end{Corollary}

From now on, $\op{char}(k)=2$ and $B$ is a real $2$-block that has a defect group $D\cong D_{2^d}$. R. Brauer \cite{BrauerDihedral} showed that $B$ has $2^{d-2}+3$ irreducible characters. Four of these, denoted $\chi_1,\chi_2,\chi_3,\chi_4$, have height $0$. The remaining $2^{d-2}-1$ irreducible characters have height $1$, and fall into $d-2$ disjoint families $F_0,\ldots,F_{d-3}$. The family $F_i$ consists of\/ $2^i$ irreducible characters, all of whom are $2$-conjugate (i.e. conjugate via Galois automorphisms that fix all $2$-power roots of unity).

Choose $\chi^{(j)}\in F_j$, and set $\epsilon^{(j)}=\epsilon(\chi^{(j)})$, for $j=0,\ldots,d-3$. Note that $\epsilon^{(j)}$ does not depend on the choice of\/ $\chi^{(j)}$ in $F_j$, as Galois conjugation preserves FS-indicators. Also set $\epsilon_i=\epsilon(\chi_i)$, for $i=1,2,3,4$. 

We also fix a defect couple $(D,E)$ of\/ $B$. If\/ $B$ is principal, then $E=D$ and we abuse notation by saying that $E$ is of type (a). Otherwise $[E:D]=2$ and $E$ is of type (a),(b),(c),(d) or (e), as in Proposition \ref{P:Dihedral_extensions}.

We give representatives $x$ for the $E$ conjugacy classes in $E\backslash D$ in the table below. We also indicate when $x$ is an involution:
\begin{table}[h]\begin{center}
\begin{tabular}{|c|c|c|c|}
\hline
$E$ & $x\in E\backslash D$ & $o(x)=2$ & $\op{C}_D(x)$\\ \hline
 (a) or (b) & $e$ & (a) & $D$\\
    & $s_1e$ & (a) & $D$\\
    & $s^ie,\, \scriptstyle{1\leq i\leq2^{d-2}\!-1}$ & $s_2e$ if (b) & $S$\\
    & $te$ & (a) & $X_2$\\
    & $ste$ & (a) & $Y_2$\\ \hline
 (c) or (d)   & $s^ite,\, \scriptstyle{0\leq i\leq2^{d-2}\!-1}$ & never & $S$\\
    & $e$ & (c) & $S_1$\\ \hline
 (e)    & $e$ & (e) & $X_{d-1}$\\
    & $s_2e$ & never & $Y_{d-1}$\\
    & $s^ie,\, \scriptstyle{1\leq i\leq2^{d-3}\!-1}$ & never & $S_{d-2}$\\
    & $te$ & (e) & $X_2$\\
    & $ss_2te$ & never &$Y_2$\\ \hline
\end{tabular}
\end{center}\caption{}\label{Tb:E}\end{table}


\section{Local structure and Real Characters}\label{S:local}

We first enumerate the real subpairs and the real irreducible characters of height $1$, in a $2$-block with dihedral defect group. Fix a sylow $B$-subpair $(D,b_D)$, such that $E\leq\op{N}^*(b_D)$. So $(D,E)$ is a defect couple of\/ $b_D^{E\op{C}(D)}$.

For each $Q\leq D$, there is a unique subpair $(Q,b_Q)\leq(D,b_D)$. The notation $b_Q$ is compatible with that used in \cite{BrauerDihedral}. We use $\sim$ or $\sim_G$ to denote $G$-conjugacy of subpairs.

\begin{Lemma}\label{L:subpair_conjugacy}
Let $Q\leq D$, $Q\not\leq S$ with $|Q|=2^i$. If\/ $Q\leq X_{d-1}$, then $(Q,b_Q)\sim_D(X_i,b_{X_i})$. If\/ $Q\leq Y_{d-1}$, then $(Q,b_Q)\sim_D(Y_i,b_{Y_i})$. Along with these, the following generate all $G$-conjugacies among the $B$-subpairs contained in $(D,b_D)$:
\begin{itemize}
\item[(i)] if  $N(b_{X_2})/\op{C}(X_2)\cong {\mathfrak S}_3$ then $(\langle t\rangle,b_{\langle t\rangle})\sim(S_1,b_{S_1})$.
\item[(ii)] if  $N(b_{Y_2})/\op{C}(Y_2)\cong {\mathfrak S}_3$ then $(\langle st\rangle,b_{\langle st\rangle}\!)\sim(S_1,b_{S_1})$.
\end{itemize}
\end{Lemma}

Following Brauer, case (aa) is the simultaneous occurrence of (i) and (ii); case (ab) is the occurrence of (i) but not (ii); case (ba) is the occurrence of (ii) but not (i); case (bb) is the occurrence of neither (i) nor (ii). Then $B$ has $3$, $2$ or $1$ irreducible modules according as case (aa), (ab) (or (ba)) or (bb) occurs. Moreover, $B$ is a nilpotent block in case (bb).

\begin{Lemma}\label{L:cd13}
If\/ $E$ is of type (c) or (d) then $l(B)\ne 2$.
\begin{proof}
There exists $f\in E\backslash D$ such that $X_2^f=Y_2$. So $\op{C}(X_2)^f=\op{C}(Y_2)$. But
$(X_2,b_{X_2})^f=(Y_2,b_{Y_2}^o)$, by Lemma \ref{L:eo_action}. Thus $\op{N}(b_{X_2})^f=\op{N}(b_{Y_2}^o)=\op{N}(b_{Y_2})$. The conclusion now follows, as $N(b_{X_2})/\op{C}(X_2)\cong N(b_{Y_2})/\op{C}(Y_2)$.
\end{proof}
\end{Lemma}

The following three lemmas can be proved by careful applications of Theorem \ref{T:real_subpairs}, using the information in Table \ref{Tb:E}. We omit the proofs.

\begin{Lemma}\label{L:real_Sylow_subpair}
$(D,b_D)$ is real if and only if\/ $E$ is of type (a) or (b).
\end{Lemma}

\begin{Lemma}\label{L:real_S_1}
$(S_1,b_{S_1})$ is real. Moreover, it is strongly real if and only if\/ $E$ is not of type (d).
\end{Lemma}

\begin{Lemma}\label{L:real_subpairs}
Suppose that $1<Q<D$. Then $(Q,b_Q)$ is real if and only if one of the following holds:
\begin{itemize}
\item $E$ is of type (a) or (b).
\item $E$ is of type (c) or (d) and $(Q,b_Q)\sim(S_i,b_{Y_i})$ for some $i$.
\item $E$ is of type (e) and $(Q,b_Q)\not=(S,b_S)$.
\end{itemize}
Of these, $(Q,b_Q)$ is strongly real in the following cases:
\begin{itemize}
\item $E$ is of type (a).
\item $E$ is of type (b) and $(Q,b_Q)\sim(S_i,b_{S_i})$, for some $i$.
\item $E$ is of type (c) and $(Q,b_Q)\sim(S_1,b_{S_1})$.
\item $E$ is of type (e) and $(Q,b_Q)\sim(S_i,b_{S_i})$ or $(X_i,b_{X_i})$ for some $i$.
\end{itemize}
\end{Lemma}

For $d\in D$ we set $b_d:=b_{\langle d\rangle}$. So $(\langle d\rangle,b_d)\leq(D,b_D)$ is a $B$-subsection. According to \cite{BrauerDihedral}, if\/ $d\ne 1$ then $b_d$ has a unique irreducible Brauer character $\theta$. We use the notation $d_{\chi\theta}^{(d)}$ for the generalized decomposition number associated with $\chi\in\op{Irr}(B)$ i.e. we supress the dependence of\/ $\theta$ on $d$. The associated $B$-column is $(d,\theta)$.

\begin{Lemma}\label{L:subsection_conjugacy}
The following generate all $D$-conjugacies among\newline $B$-subsections:
\begin{itemize}
\item $(d,b_d)\sim(d^{-1},b_d)$, for all $d\in D$.
\item if\/ $d\in X_{n-1}\backslash S_1$ then $(d,b_d)\sim(t,b_t)$.
\item if\/ $d\in Y_{n-1}\backslash S_1$ then $(d,b_d)\sim(st,b_{st})$.
\end{itemize}
Together with $D$-conjugacies, the following generate all $G$-conjugacies among $B$-subsections:
\begin{itemize}
\item[(i)]  if  $N(b_{X_2})/\op{C}(X_2)\cong {\mathfrak S}_3$ then $(t,b_t)\sim(s_1,b_{s_1})$.
\item[(ii)] if  $N(b_{Y_2})/\op{C}(Y_2)\cong {\mathfrak S}_3$ then $(st,b_{st})\sim(s_1,b_{s_1})$.
\end{itemize}
\begin{proof}
Suppose that $d\in S\backslash S_1$. Then $D\leq\op{N}(b_x)$ and $d^t=d^{-1}$. So $(d,b_d)^t=(d^{-1},b_d)$. The first statement follows from this. This was already proved by Brauer in \cite[(4.16)]{BrauerDihedral}. The next two statements follow from the discussion before Lemma \ref{L:subpair_conjugacy}. The remaining statements follow from Proposition (4A) of \cite{BrauerDihedral}.
\end{proof}
\end{Lemma}

Brauer showed in \cite{BrauerDihedral} that $B$ has $2^{d-2}+3$ columns. Exactly $5$ of these are $2$-rational, namely those of the form $(d,\theta)$, with $d=1,s_1,s_2,t$ or $st$ and $\theta\in\op{IBr}(\op{C}(d))$. There are $d-2$ families, with representatives $\{(s_i,\theta)\mid i=2,\ldots,d-1\}$. The family of\/ $(s_i,\theta)$ contains the $2^{i-2}$ columns $\{(s_i^r,\theta)\mid 1\leq r\leq 2^{i-1}-1,\mbox{ $r$ odd}\}$, which form a single $2$-conjugate orbit. Recall the notation $\chi_i,\chi^{(i)}$ for the irreducible characters in $B$. The $5$ irreducible characters which are $2$-rational are $\chi_1,\chi_2,\chi_3,\chi_4$, and $\chi^{(0)}$.

We can now prove the main result of this section.

\begin{Theorem}\label{T:real_characters}
The number of real $2$-rational irreducible characters in $B$ equals the number of real $2$-rational columns in $B$. For $d\geq4$, all irreducible characters in $F_0,\ldots,F_{d-4}$ are real, while all irreducible characters in $F_{d-3}$ are real if and only if\/ $E$ is not of type (e).
\begin{proof}
The first statement is trivial if\/ $d=3$. So we assume that $d\geq4$.

Let $m$ be the largest odd divisor of\/ $|G|$ and let $\omega$ be a primitive $(2^{d-1}m)$-th root of unity. Then $\omega=\omega_2\omega_{2'}$, where $\omega_2$ is a primitive $2^{d-1}$-th root of unity and $\omega_{2'}$ is a primitive $m$-th root of unity. Let $\gamma,\sigma\in\op{Gal}({\mathbb Q}(\omega)/{\mathbb Q})$ be such that $\omega^\gamma=\omega_2^5\omega_{2'}$ and $\omega^\sigma=\overline\omega$. Then $\gamma$ has order $2^{d-3}$ and $\sigma$ is an involution. Set $\tau:=\gamma^{2^{d-4}}$. So $\tau$ is the unique involution in $\langle\gamma\rangle$.

Set ${\mathcal G} :=\langle\gamma\rangle\times\langle\sigma\rangle$ and let $\alpha\in{\mathcal G}$. There is an integer $(\alpha)$ such that $\omega^\alpha=\omega^{(\alpha)}$. For $g\in G$, set $g^\alpha:=g^{(\alpha)}$. If\/ $\theta$ is a character (ordinary or Brauer) of a subgroup of\/ $G$, define $\theta^\alpha$ by $\theta^\alpha(g):=\theta(g^{\alpha})$, for all $g$ in the domain of\/ $\theta$. Then set $(d,\theta)^{\alpha}:=(d^{\alpha},\theta^{\alpha})$. In this way ${\mathcal G}$ acts on the ordinary irreducible characters and on the columns of\/ $B$. Note that a character is $2$-rational if it is fixed by $\gamma$, and real if it is fixed by $\sigma$.

It can be checked that $d_{\chi^{\alpha}\theta^\alpha}^{(x^{\alpha})}=d_{\chi\theta}^{(x)}$, for all $\alpha\in{\mathcal G}$. As the generalized decomposition matrix is non-singular, it follows that the characters of\/ ${\mathcal G}$ acting on the irreducible characters and on the columns coincide. We cannot use Brauer's permutation lemma to deduce that the permutation actions are isomorphic, as ${\mathcal G}$ is not cyclic. However, we can prove this fact using a careful case-by-case analysis.

Let $i=1,\ldots,d-1$. Lemmas \ref{L:real_subpairs} and \ref{L:subsection_conjugacy} imply that $(s_i,\theta)$ is not real if and only if\/ $E$ is of type (e) and $i=d-1$. In particular the $2$-rational columns $(s_1,\theta)$ and $(s_2,\theta)$ are necessarily real. Moreover, $B$ has at least one real irreducible Brauer character (this holds for any real block).

We distinguish $5$ cases:

Case (I) is that $E$ is not of type (e) and all $2$-rational columns are real. Then all columns are real. We conclude from Theorem \ref{T:real_characters} that all irreducible characters in $B$ are real.

Case (II) is that $E$ is not of type (e) and two $2$-rational columns are nonreal. Then there are two nonreal columns. We deduce from Theorem \ref{T:real_characters} that there are two nonreal irreducible characters. Consider the action of\/ $\langle\gamma\sigma\rangle$ on columns and on irreducible characters. The stabilizer of a nonreal column is $\langle\gamma^2\rangle$. The remaining $2$-rational columns are fixed by $\gamma\sigma$. For $i=3,\ldots,d-1$, the stabilizer of\/ $(s_i,\theta)$ is $\langle\gamma^{2^{i-2}}\rangle$. We deduce that exactly two orbits on the columns have stabilizer $\langle\gamma^2\rangle$. It then follows from Brauer's permutation lemma that two orbits on the irreducible characters have stabilizer $\langle\gamma^2\rangle$. A $2$-rational orbit has stabilizer $\langle\gamma^2\rangle$ if and only if it is nonreal. There is at most one such orbit. Suppose that a family $F_i$ has stabilizer $\langle\gamma^2\rangle$, for $i\geq2$. Then $|F_i|\leq[\langle\gamma\rangle:\langle\gamma^2\rangle]=2$. So $i=2$. We conclude that there is a pair of nonreal $2$-rational characters in $B$. Moreover, $F_1$ has stabilizer $\langle\gamma^2\rangle$, whence it is real.

Case (III) is that $E$ is of type (e) and all $2$-rational columns are real. Then there are exactly $2^{d-3}$ nonreal columns, namely the $2$-conjugates of\/ $(s,\theta)$. We deduce from Theorem \ref{T:real_characters} that there are $2^{d-3}$ nonreal irreducible characters. As $\langle\sigma\rangle$ acts nontrivially on the $2$-conjugates of\/ $(s,\theta)$, this column has stabilizer $\langle\tau\sigma\rangle$ in ${\mathcal G}$. All other columns are fixed by $\tau\sigma$. So by Brauer's permutation lemma all irreducible characters are $\langle\sigma\tau\rangle$-invariant. Let $\chi\in F_{d-3}$. Then $\chi^{\sigma}=\chi^{\tau}\ne\chi$. So  $F_{d-3}$ is a nonreal family of characters. This accounts for all nonreal irreducible characters in $B$.

Case (IV) is that $d=4$, $E$ is of type (e) and exactly two $2$-rational columns are nonreal. Thus $B$ has two nonreal $2$-rational columns and $2$ nonreal $2$-irrational columns. It then follows from Theorem \ref{T:real_characters} that $B$ has $4$ nonreal irreducible characters. As $B$ has $4$ irreducible characters of height $0$ and $3$ of height $1$, it follows that $B$ has a pair of nonreal characters of height $0$ and another pair of height $1$. The latter belong to $F_1$.

Case (V) is that $d>4$, $E$ is of type (e) and exactly two $2$-rational columns are nonreal. So $B$ has $2^{d-3}+2$ nonreal columns and hence $2^{d-3}+2$ nonreal irreducible characters. We consider the action of\/ $\langle\gamma\sigma\rangle$ on columns and on irreducible characters. Just as in case (II), there are two orbits on the columns whose stabilizer is $\langle\gamma^2\rangle$ (as $d>4$, the stabilizer $\langle\tau\sigma\rangle$ of\/ $(s_{d-1},\theta)$ in $\langle\gamma\sigma\rangle$ has trivial intersection with $\langle\gamma\sigma\rangle$). Again just as in case (II), we conclude that there is a pair of nonreal $2$-rational characters in $B$. Suppose that some character in $F_i$ has a trivial stabilizer. Then $|F_i|\geq|\langle\gamma\sigma\rangle|=2^{d-3}$. So $F_i=F_{d-3}$. Let $\chi\in F_i$. Then $\chi^{\sigma}=\chi^{\gamma^{-1}}\ne\chi$. So $F_{d-3}$ contains $2^{d-3}$ nonreal characters. As we have accounted for all $2^{d-3}+2$ nonreal irreducible characters in $B$, the proof is complete.
\end{proof}
\end{Theorem}

\section{Morita equivalence classes}\label{S:Morita}

Karin Erdmann classified the possible blocks with dihedral defect group, up to morita equivalence, in \cite{ErdmannTame}. Her results are summarized in the table on pp294--296. It is not known if each of the possible $8$ Morita equivalence classes occur for all dihedral groups.

Let $B$ be a $2$-block with a dihedral defect group $D$, where $|D|=2^d$. Then by \cite{BrauerDihedral} $k(B)=2^{d-2}+3$ and $l(B)\leq3$. If $l(B)=1$ then $B$ is a nilpotent block. In particular $B$ is Morita equivalent to $kD$. This corresponds to Brauer's case (bb). The case $l(B)=3$ corresponds to Brauer's case (aa). Then there are three possible Morita equivalence classes. The principal $2$-blocks of the groups $\op{L}_2(q)$ with $q\equiv 1$ or $3$ mod $4$, and the principal block of the alternating group ${\mathfrak A}_7$ provide examples of each class. As Erdmann remarks on p293 of \cite{ErdmannTame}, her proof in \cite{ErdmannDihedral} that $d=3$ for the last class is erroneous. On the other hand, there are no known examples (at least to this author) of such blocks with $d>3$.

Finally we consider when $l(B)=2$, which corresponds to Brauer's cases (ab) and (ba). Then there are quivers $Q_1,Q_2$, a parameter $c=0,1$, and ideals $I_i(c)\leq Q_i$, $i=1,2$ inside the ideal of paths of length $\geq2$, such that $B$ is Morita equivalent to exactly one of the $4$ algebras $kQ_i/I_i(c)$. By \cite{ErdmannGU}, the case $kQ_2/I(0)$ occurs for the principal blocks of certain quotients of the unitary groups $\op{GU}_2(q)$, when $q\equiv3$ mod $4$. By a new result of \cite{HolmZhou}, both $kQ_1/I_1(1)$ and $kQ_2/I_2(1)$ occur for the principal blocks of $\op{PGL}_2(q)$, corresponding to $q\equiv1$ or $-1$ mod $8$, respectively. Note that in these cases $d\geq4$. It turns out that in all cases the decomposition matrix of $B$ depends on $Q_i$, but not on $c$. For this reason, we will refer to $B$ as being of type $\op{PGL}_2(q)$, with $q\equiv1$ mod $4$, for $Q_1$, or $q\equiv3$ mod $4$, for $Q_2$.

The papers \cite{BrauerDihedral}, \cite{Butler}, \cite{Donovan}, \cite{DonovanPreislich} and \cite{ErdmannDihedral} give more information on the decomposition matrices and modules of these blocks. Derived equivalences between the block algebras are discussed in \cite{Holm} and \cite{Linckelmann}.

We obtain the form of the generalized decomposition matrices given in the next lemma using \cite{BrauerDihedral}, \cite{Donovan} and \cite{ErdmannDihedral}. In particular, we make use of the parameters $\delta_1,\delta_2,\delta_3,\delta_4$ introduced in Proposition (6C) of \cite{BrauerDihedral}, and the further information given by Theorem 5 and Propositions (6H) and (6I) of that paper. The $\delta_i$ are evaluated when $l(B)=2$ in \cite{Donovan} (implicitly) and when $l(B)=3$ in \cite{ErdmannDihedral}. Entries that we don't care to specify are denoted by $*$.

\begin{Lemma}\label{L:MoritaTypes}
For $d\geq3$, a block $B$ with a dihedral defect group of order $2^d$ has one of $6$ possible generalized decomposition matrices. Allowing $i$ to range over $0,\ldots,d-4$, and $j$ over $2,\ldots,d-2$, there are signs $\varepsilon,\varepsilon',\varepsilon_j$ such that this matrix has the form:
\begin{enumerate}
\item[(i)] $B$ is nilpotent and has one irreducible module $M_1$.
$$
\begin{array}{l|rrrrrr}
            & M_1 & s_1 & t & st & s_j & s\\
\hline
\chi_1      & 1 &  \varepsilon & 1 & \varepsilon' & \varepsilon_j & 1\\
\chi_2      & 1 &  \varepsilon &-1 &-\varepsilon' & \varepsilon_j & 1\\
\chi_3      & 1 &  \varepsilon & 1 &-\varepsilon' & \varepsilon_j &-1\\
\chi_4      & 1 &  \varepsilon &-1 & \varepsilon' & \varepsilon_j &-1\\
\chi^{(i)}  & 2 & 2\varepsilon & 0 & 0 		   & *             & *\\
\chi^{(d-3)}& 2 &-2\varepsilon & 0 & 0 		   & *             & *
\end{array}
$$
\item[(ii)] $B$ has two irreducible modules $M_1,M_2$ and is of type $\op{PGL}(2,q)$ with $q\equiv 1$ mod $4$.
$$
\begin{array}{l|rrrrrr}
            & M_1 & M_2 & s_1 & t & s_j & s\\
\hline
\chi_1      & 1 & 0 &  \varepsilon & 1 & \varepsilon_j & 1 \\
\chi_2      & 1 & 1 &  \varepsilon &-1 & \varepsilon_j & 1\\
\chi_3      & 1 & 0 &  \varepsilon &-1 & \varepsilon_j &-1\\
\chi_4      & 1 & 1 &  \varepsilon & 1 & \varepsilon_j &-1\\
\chi^{(i)}  & 2 & 1 & 2\varepsilon & 0 & *             & *\\
\chi^{(d-3)}& 2 & 1 &-2\varepsilon & 0 & *             & *
\end{array}
$$
\item[(iii)] $B$ has two irreducible modules $M_1,M_2$ and is of type $\op{PGL}(2,q)$ with $q\equiv 3$ mod $4$.
$$
\begin{array}{l|rrrrrr}
            & M_1 & M_2 &  s_1 & t & s_j & s\\
\hline
\chi_1      & 1 & 0 & -\varepsilon & 1 & -\varepsilon_j &-1\\
\chi_2      & 1 & 1 &  \varepsilon &-1 &  \varepsilon_j & 1\\
\chi_3      & 1 & 0 & -\varepsilon &-1 & -\varepsilon_j & 1\\
\chi_4      & 1 & 1 &  \varepsilon & 1 &  \varepsilon_j &-1\\
\chi^{(i)}  & 0 & 1 & 2\varepsilon & 0 & *             & *\\
\chi^{(d-3)}& 0 & 1 &-2\varepsilon & 0 & *             & *
\end{array}
$$
\item[(iv)] $B$ has $3$ irreducible modules $M_1,M_2,M_3$. If\/ $d=3$, then $B$ is Morita equivalent to the principal block of\/ $A_7$.
$$
\begin{array}{l|rrrrrr}
            & M_1 & M_2 & M_3 & s_1 & s_j & s\\
\hline
\chi_1      & 1 & 0 & 0 & -\varepsilon & -\varepsilon_j &-1\\
\chi_2      & 1 & 1 & 0 &  \varepsilon &  \varepsilon_j & 1\\
\chi_3      & 1 & 1 & 1 &  \varepsilon &  \varepsilon_j &-1\\
\chi_4      & 1 & 0 & 1 & -\varepsilon & -\varepsilon_j & 1\\
\chi^{(i)}  & 0 & 1 & 0 & 2\varepsilon & *             & *\\
\chi^{(d-3)}& 0 & 1 & 0 &-2\varepsilon & *             & *
\end{array}
$$
\item[(v)] $B$ has $3$ irreducible modules $M_1,M_2,M_3$ and is Morita equivalent to the principal block of\/ $\op{PSL}(2,q)$, with $q\equiv 1$ mod $4$.
$$
\begin{array}{l|rrrrrr}
            & M_1 & M_2 & M_3 & s_1 & s_j & s\\
\hline
\chi_1      & 1 & 0 & 0 &  \varepsilon &  \varepsilon_j & 1\\
\chi_2      & 1 & 1 & 1 &  \varepsilon &  \varepsilon_j & 1\\
\chi_3      & 1 & 1 & 0 &  \varepsilon &  \varepsilon_j &-1\\
\chi_4      & 1 & 0 & 1 &  \varepsilon &  \varepsilon_j &-1\\
\chi^{(i)}  & 2 & 1 & 1 & 2\varepsilon & *             & *\\
\chi^{(d-3)}& 2 & 1 & 1 &-2\varepsilon & *             & *
\end{array}
$$
\item[(vi)] $B$ has $3$ irreducible modules $M_1,M_2,M_3$ and is Morita equivalent to the principal block of\/ $\op{PSL}(2,q)$, with $q\equiv 3$ mod $4$.
$$
\begin{array}{l|rrrrrr}
            & M_1 & M_2 & M_3 & s_1 & s_j & s\\
\hline
\chi_1      & 1 & 0 & 0 & -\varepsilon & -\varepsilon_j & -1\\
\chi_2      & 1 & 1 & 1 &  \varepsilon &  \varepsilon_j & 1\\
\chi_3      & 0 & 1 & 0 &  \varepsilon &  \varepsilon_j & -1\\
\chi_4      & 0 & 0 & 1 &  \varepsilon &  \varepsilon_j & -1\\
\chi^{(i)}  & 0 & 1 & 1 & 2\varepsilon & *              & *\\
\chi^{(d-3)}& 0 & 1 & 1 &-2\varepsilon & *              & *
\end{array}
$$
\end{enumerate}
\begin{proof}
Brauer showed in \cite{BrauerDihedral} that $d^{(s_1)}_{\chi,\theta}=\pm 1$ or $\pm 2$, depending on whether $\chi$ has height $0$ or $1$. Moreover, $d^{(s_1)}_{\chi,\theta}$ is constant on each of the families of\/ $2$-conjugate characters $F_1,\ldots,F_{d-3}$. If\/ $B$ is of type (bb) or (ab), and $x=t,st$, then $d^{(x)}_{\chi,\theta}=\pm 1$ or $0$, depending on whether $\chi$ has height $0$ or $1$. The columns for $s_1,t,st,s_j,s$ can then be recovered, up to the signs $\varepsilon,\varepsilon',\varepsilon_j$, and possible rearranging of the $\chi_i$, using orthogonality with the columns of the decomposition matrix of\/ $B$.
\end{proof}
\end{Lemma}

\begin{Corollary}\label{C:height1FS-1}
Suppose that $B$ has a dihedral defect group of order at least $16$ and that all height $1$ irreducible characters in $B$ are real-valued. Then at least one of these characters has FS-indicator $+1$.
\begin{proof}
The multiplicity of an irreducible $B$-module $M$ in $k\invol B$ is given by $\sum_{\chi\in\op{Irr}(B)}\epsilon(\chi)d_{\chi,M}$. In particular this sum is non-negative. The hypothesis on the defect group implies that $B$ has at least $3$ irreducible characters of height $1$. We assume that all irreducible characters in $B$ of height $1$ have FS-indicator $-1$, and derive a contradiction.

Suppose that $B$ is nilpotent. So $M=M_1$. The contribution of the height zero characters to the above sum is at most $+4$. On the other hand, the height $1$ characters contribute an integer $\leq-6$, whence the sum is negative.

Suppose that $B$ is not nilpotent. Take $M=M_2$. The contribution of the height zero characters to the above sum is at most $+2$. On the other hand, the height $1$ characters contribute an integer $\leq-3$, whence the sum is negative.
\end{proof}
\end{Corollary}

We will use the following rather technical result:

\begin{Lemma}\label{L:s^2}
Let $f$ be a ${\mathbb C}$-valued function on $\op{Irr}(B)$ that is constant on each of $F_1,\ldots,F_{d-3}$. Then
$$
\sum_{\chi\in\op{Irr}(B)}f(\chi)d_{\chi,\theta}^{(s^2)}=
\sum_{i=1}^4f(\chi_i)d_{\chi_i,\theta}^{(s^2)}+f(\chi^{(0)})d_{\chi^{(0)},\theta}^{(s^2)}.
$$
\begin{proof}
Let $1\leq r\leq d-3$. Then $\zeta+\zeta^{-1}$ takes on $2^{r-1}$ different values as $\zeta$ ranges over the primitive $2^{r+1}$-th roots of unity. A careful reading of Section 6 of \cite{BrauerDihedral} shows that for each $\zeta$ there are exactly two irreducible characters $\chi\in F_r$ with $d_{\chi,\theta}^{(s^2)}=\zeta+\zeta^{-1}$. Now $f(\chi)$ is constant on $F_r$. So the net contribution of the two characters associated with each of\/ $\zeta$ and $-\zeta$ to $\sum_{\chi\in\op{Irr}(B)}f(\chi)d_{\chi,\theta}^{(s^2)}$ is $0$. The lemma follows from this.
\end{proof}
\end{Lemma}

\section{Type {\rm (a)}: totally split extended defect groups}\label{S:TotallySplit}

In this section $E=D\times\langle e\rangle$ i.e. $e\in\op{C}(D)$ and $e^2=1_G$. This includes the principal $2$-block case, when $e=1_G$. By Theorem \ref{T:real_characters}, all height $1$ irreducible characters in $B$ are real-valued. Using Lemma \ref{L:real_subpairs} and Theorem \ref{T:real_characters}, we can show that if\/ $B$ is nilpotent then all its irreducible characters are real valued. This also follows from Lemma 2.2 of \cite{Murray_GroupT}.

\begin{Lemma}\label{L:TotallySplitVertexD}
There is an indecomposable $B$-module $M_D$ such that $M_D$ has vertex $D$ and $M_D\oplus M_D$ is the sum of all components of\/ $k\invol B$ that have vertex $D$. Moreover $\nu(\op{dim}(M_D))=\nu[G:D]$.
\begin{proof}
Let $b$ be the $2$-block of\/ $\op{N}(D)$ that is the Brauer correspondent of\/ $B$. Then $b$ is real and has defect couple $(D,E)$. Now $\op{N}(D)$ acts trivially on $\mathord{Z}(D)$. So by Lemma \ref{L:Alperin}, there is a self-dual irreducible $b$-module $I$, such that the sum of all components of\/ $k\invol_{\op{N}(D)}b$ that have vertex $D$ is isomorphic to $I\oplus I$. Set $M_D$ as the Green correspondent of\/ $I$, with respect to $(G,D,\op{N}(D))$. Then the sum of all components of\/ $k\invol B$ that have vertex $D$ is isomorphic to $M_D\oplus M_D$.

As $I$ is irreducible and projective as $\op{N}(D)/D$-module, $\nu(\op{dim}(I))=\nu[\op{N}(D):D]$. Now $I{\uparrow^G}=M_D\oplus W$, where every component of\/ $W$ has a vertex that is a proper subgroup of\/ $D$. In particular $\nu(\op{dim}(W))>\nu[G:D]$. We conclude that $\nu(\op{dim}(M_D))=\nu[G:D]$.
\end{proof}
\end{Lemma}

\begin{Lemma}\label{L:TotallySplitVertexXorY}
There are indecomposable $B$-modules $M_{X_2},M_{Y_2}$ such that $M_{X_2}\oplus M_{Y_2}$ is the sum of all components of\/ $k\invol B$ that have vertex $X_2$ or $Y_2$. Here $M_V$ has vertex $B$-subpair $(V,b_V)$, for $V=X_2,Y_2$.
\begin{proof}
Recall that $X_2$ is a Klein-four subgroup of\/ $D$. Set $N:=\op{N}(X_2)$ and $C:=\op{C}(X_2)$. Then $b_{X_2}$ is real, with defect couple $(X_2,X_2\times\langle e\rangle)$. As $b_{X_2}$ is nilpotent, Theorem 1.7 of \cite{Murray_GroupT} implies that $k\invol_{C}b_{X_2}\cong R^4$, where $R$ is the unique irreducible $b_{X_2}$-module. Let $\op{N}(R)$ denote the common inertia group of\/ $b_{X_2}$ and $R$ in $N$.

Set $b:=b_{X_2}^{N}$. Then $b$ is real, with defect couple $(X_3,X_3\times\langle e\rangle)$. By Lemma \ref{L:TotallySplitVertexD}, there is an indecomposable $b$-module $M_{X_3}$ such that $M_{X_3}\oplus M_{X_3}$ is the sum of all components of\/ $k\invol_{N}b$ that have vertex $X_3$. We note that $k\invol_{C}b$ is the sum of all components of\/ $k\invol_{N}b$ that have vertex $X_2$ or $X_3$.

Suppose first that $\op{N}(R)/C\cong{\mathbb Z}_2$. Then $R$ has a unique extension to $\op{N}(R)$ and $b$ is nilpotent. Let $I$ be the unique irreducible $b$-module. Then $R$ occurs once in the semisimple module $I{\downarrow_{C}}$. As $k\invol_{C}b_{X_2}\cong R^4$, it follows that $[k\invol_{C}:I]=4$. Now $\nu(\op{dim}(M_{X_3}))=\nu[N:X_3]=\nu(\op{dim}(I))$. So $[M_{X_3}:I]$ is odd. The only possibility is that $M_{X_3}\cong I$. We may now write $k\invol_{C}b\cong I\oplus I\oplus W$, where $[W:I]=2$, and every component of\/ $W$ has vertex $X_2$. As $I$ does not itself have vertex $X_2$, it follows that $W$ is indecomposable.

Suppose then that $\op{N}(R)/C\cong{\mathfrak S}_3$. Then two irreducible $\op{N}(R)$-modules lie over $R$. So $b$ has two irreducible modules $I_1$ and $I_2$. We may assume that $I_1$ has vertex $X_3$ and $I_1{\downarrow_{C}}\cong R$, while $I_2$ has vertex $X_2$ and $I_2{\downarrow_{C}}\cong R\oplus R$. By considering the Clifford theory of the ordinary characters of\/ $b$, we see that $b$ is Morita equivalence to the principal $2$-block of\/ $\op{PGL}(2,3)\cong{\mathfrak S}_4$. Arguing as in the previous paragraph, $M_{X_3}\cong I_1$ and $k\invol_{C}b\cong I_1\oplus I_1\oplus W$, where $W$ is indecomposable with vertex ${X_2}$. With more work, we can even show that $W\cong I_2$.

Regardless of the structure of\/ $\op{N}(R)/C$, let $M_{X_2}$ denote the Green correspondent of\/ $W$. Then $M_{X_2}$ is the unique component of\/ $k\invol B$ that has vertex $B$-subpair $(X_2,b_{X_2})$. Analogous arguments show that $k\invol B$ has a unique component $M_{Y_2}$ that has vertex $B$-subpair $(Y_2,b_{Y_2})$. But Lemma \ref{L:subpair_conjugacy} shows that $(X_2,b_{X_2})\not\sim(Y_2,b_{Y_2})$. So $M_{X_2}\not\cong M_{Y_2}$.
\end{proof}
\end{Lemma}

\begin{Theorem}\label{T:TotallySplitModules}
$k\invol B\cong M_D\oplus M_D\oplus M_{X_2}\oplus M_{Y_2}$. If\/ $B$ is nilpotent then all its irreducible characters have FS-indicator $+1$, $M_D$ is irreducible and $M_{X_2},M_{Y_2}$ are each of composition length $2^{d-2}$.
\begin{proof}
Suppose first that $B$ is nilpotent, with unique irreducible module $M_1$. Let $\theta$ be the Brauer character of\/ $M_1$. From the discussion in Section \ref{S:Introduction}, $M_1$ occurs with multiplicity $\sum\epsilon(\chi)d_{\chi,\theta}$ in $k\invol B$. By Part (i) of Lemma \ref{L:MoritaTypes}, this is at most $4+2\sum_{j=0}^{d-3}2^j=|D|/2+2$. Now $[M_D:M_1]\geq1$ and by consideration of vertices, $|D|/4$ divides both $[M_{X_2}:M_1]$ and $[M_{Y_2}:M_1]$. It then follows from Lemmas \ref{L:TotallySplitVertexD} and \ref{L:TotallySplitVertexXorY} that $[k\invol B:M_1]\geq|D|/2+2$. All statements of the theorem now follows for nilpotent $B$.

Now let $B$ be of arbitary Morita equivalence type and let $(V,b_V)$ be a vertex $B$-subpair of a component of\/ $k\invol B$. Lemmas \ref{L:detectVertex} and \ref{L:subpair_conjugacy} imply that $(V,b_V)$ is a vertex $b_x$-subpair of a component of\/ $k\invol_{\op{C}(x)}b_x$, where $x=s_1,t$ or $st$. Now $b_{s_1}$ is real and nilpotent, with defect couple $(D,D\times\langle e\rangle)$. So by the previous paragraph,  $(V,b_V)$ is conjugate to $(D,b_D),(X_2,b_{X_2})$ or $(Y_2,b_{Y_2})$, if\/ $x=s_1$. The block $b_t$ is real with defect couple $(X_2,X_2\times\langle e\rangle)$ and Sylow $B$-subpair $(X_2,b_{X_2})$. As $X_2$ is a Klein-four group, Theorem 1.7(i) of \cite{Murray_GroupT} implies that $(V,b_V)$ is conjugate to $(X_2,b_{X_2})$, if\/ $x=t$. In the same way, $(V,b_V)$ is conjugate to $(Y_2,b_{Y_2})$, if\/ $x=st$. The first statement of the Theorem now follows from Lemmas \ref{L:TotallySplitVertexD} and \ref{L:TotallySplitVertexXorY}.
\end{proof}
\end{Theorem}

\begin{Theorem}\label{T:TotallySplit}
Let $\chi\in\op{Irr}(B)$. Then $\nu(\chi)=+1$, unless $B$ is of type (vi) and $\chi=\chi_3$ or $\chi_4$; in that case $\ov{\chi_3}=\chi_4$.
\begin{proof}
As $b_{s_1}$ is nilpotent, real, and has defect couple $(D,D\times\langle e\rangle)$, Theorem \ref{T:TotallySplitModules} shows that $\epsilon(\psi)=+1$, for all $\psi\in\op{Irr}(b_{s_1})$. From our knowledge of\/ $d_{\psi,\theta}^{(s_1)}$, and the positivity assertion in Lemma \ref{L:BrauerSum}, we get
\begin{equation}\label{E:BrauerSumType(i)}
\sum_{\psi\in\op{Irr}(b_{s_1})}\epsilon(\psi)d_{\psi,\theta}^{(s_1)}=\pm(4+2(2^{d-3}-1)-2.2^{d-3})=2.
\end{equation}

Let $d\geq 4$. Then $b_{s^2}$ is real, with defect couple $(S,S\times\langle e\rangle)$. Theorem 1.6 of \cite{Murray_GroupT} shows that $b_{s^2}$ has two real-valued irreducible characters, $\psi_1$ and $\psi_2$. Moreover, $\epsilon(\psi_i)=+1$, for $i=1,2$. It follows that
\begin{equation}\label{E:BrauerSumType(i)s^2}
\begin{array}{lcl}
\sum\limits_{\psi\in\op{Irr}(b_{s^2})}\epsilon(\psi)d_{\psi,\theta}^{(s^2)}
	&=&\pm\epsilon(\psi_1)\pm\epsilon(\psi_2)\\
	&=&+2,\quad\mbox{by Lemma \ref{L:2ElementSum}, as $s^2=(se)^2$.}
\end{array}
\end{equation}

Suppose first that $B$ has type (ii) or (v). Lemma \ref{L:MoritaTypes} gives the form of the generalized decomposition matrix of $B$. Applying Lemma \ref{L:BrauerSum} to the column $d_{\chi\theta}^{(s_1)}$, equality \eqref{E:BrauerSumType(i)} gives
$$
2=\sum_{\chi\in\op{Irr}(B)}\epsilon(\chi)d_{\chi,\theta}^{(s_1)}=\varepsilon(\sum_{i=1}^4\epsilon_i+2\sum_{j=0}^{d-4}\epsilon^{(j)}2^j-2\epsilon^{(d-3)}2^{d-3}).
$$
In particular $\sum_{i=1}^4\epsilon_i\equiv 0$ mod $4$. It then follows from Corollary \ref{C:heightzero_FS+1} 
that $\epsilon_i=1$, for $i=1,2,3,4$. Then the above equation rearranges to
$$
\sum_{j=0}^{d-4}\epsilon^{(j)}2^j-\epsilon^{(d-3)}2^{d-3}\,=\,\varepsilon-2\in\{-1,-3\}.
$$
By Lemma \ref{L:SignedSum} the two solutions: $\varepsilon=+1$ and $\epsilon^{(j)}=+1$, for $j\geq0$, or $\varepsilon=-1$, $\epsilon^{(0)}=-1$ and $\epsilon^{(j)}=+1$ for $j>0$, are the only ones. We claim that the latter solution does not occur. For, taking $f(\chi)=d_{\chi,\theta}^{(s_1)}$ in Lemma \ref{L:s^2}, we get $\sum_{\chi\in\op{Irr}(B)}d_{\chi,\theta}^{(s_1)}d_{\chi,\theta}^{(s^2)}=4\varepsilon\varepsilon_{d-2}+2\varepsilon d_{\chi^{(0)},\theta}^{(s^2)}$. This sum is zero, by column orthogonality. So
$d_{\chi^{(0)},\theta}^{(s^2)}=-2\varepsilon_{d-2}$. Now
$$
2=\sum\limits_{\psi\in\op{Irr}(b_{s^2})}\epsilon(\psi)d_{\psi,\theta}^{(s^2)}=\sum\limits_{\chi\in\op{Irr}(B)}\epsilon(\chi)d_{\chi,\theta}^{(s^2)}=\varepsilon_{d-2}(4-2\epsilon^{(0)}),
$$
by \eqref{E:BrauerSumType(i)s^2}, Lemma \ref{L:BrauerSum}, and an application of Lemma \ref{L:s^2}, with $f(\chi)=\epsilon(\chi)$. Thus $\varepsilon_{d-2}=+1$ and $\epsilon^{(0)}=+1$, which proves the claim.

Suppose then that $B$ has type (iii), (iv) or (vi). Then from the decomposition matrices given in Lemma \ref{L:MoritaTypes}, we have $\sum_{i=1}^4\!\epsilon_id_{\chi,\theta}^{(s_1)}\!=0,\pm2$. However, just as above, this sum is $\equiv0$ modulo $4$. So it must equal $0$. If\/ $B$ is of type (iii) or (iv), we get $\epsilon_i=+1$, for all $i$; if\/ $B$ is of type (vi), we get  $\epsilon_1=\epsilon_2=+1$ and $\epsilon_3=\epsilon_4=0$. 

From knowledge of\/ $d_{\chi,\theta}^{(s_1)}$, Lemmas \ref{L:BrauerSum}, \eqref{E:BrauerSumType(i)} and the work above give
$$
2\varepsilon(\sum_{j=0}^{d-4}\epsilon^{(j)}2^j-\epsilon^{(d-3)}2^{d-3})=+2.
$$
Thus $\epsilon^{(j)}=-\varepsilon=\mp1$, for all $j=0,\ldots,d-3$.

We claim that $\varepsilon=-1$ and $\epsilon^{(j)}=+1$, for all $j$. If\/ $|D|\geq16$, this is a consequence of Corollary \ref{C:height1FS-1}. So from now on we assume that $|D|=8$ and $\epsilon^{(0)}=-1$, and argue to a contradiction.

Recall that $k\invol B=M_D^2\oplus M_{X_2}\oplus M_{Y_2}$. Moreover, $M_{X_2}\not\cong M_{Y_2}$, as $(X_2,b_{X_2})\not\sim(Y_2,b_{Y_2})$. So $M_{X_2}$ and $M_{Y_2}$ are self-dual. Also $[M_D:M_1]$ is odd. We claim that $M_D\cong M_1$. Otherwise $M_D$ is reducible, whence by Lemma \ref{L:Indecomposable}, $M_D^2$ has a composition factor that occurs with multiplicity $\geq4$. Using \eqref{E:Smult}, and Lemma \ref{L:MoritaTypes}, $[k\invol:M_1]\leq4$, ${[k\invol:M_2]=1}$ and ${[k\invol:M_3]\leq2}$. So $[M_D:M_1]=2$. This contradiction proves our claim.

If\/ $B$ is of type (iii), then ${[M_{X_2}+M_{Y_2}:M_1]=2}$ and $[M_{X_2}+M_{Y_2}:M_2]=1$, using \eqref{E:Smult} and Lemma \ref{L:MoritaTypes}. We choose notation so that $M_{X_2}$ is irreducible. So $M_{X_2}\cong M_2$. Now the ${\mathcal O}$-lift of\/ $M_{X_2}$ has character $\chi^{(0)}$. So $\chi^{(0)}$ appears with multiplicity $1$ in the permutation character of\/ ${\mathcal O}\invol$. But then $\epsilon(\chi^{(0)})=+1$, contrary to hypothesis.

If\/ $B$ is of type (iv), then ${[M_{X_2}+M_{Y_2}:M_1]=2}$, $[M_{X_2}+M_{Y_2}:M_2]=1$ and ${[M_{X_2}+M_{Y_2}:M_2]=2}$, using \eqref{E:Smult} and Lemma \ref{L:MoritaTypes}. Choose notation so that $[M_{X_2}:M_3]\ne0$. Then $\chi_3$ or $\chi_4$ occurs in the character of the ${\mathcal O}$-lift of\/ $M_{X_2}$. So $[M_{X_2}:M_1]\ne0$. As $M_{X_2}$ is reducible, Lemma \ref{L:Indecomposable} implies that $[M_{X_2}:M_1]=2$. Then $[M_{Y_2}:M_1]=0$, whence $[M_{Y_2}:M_3]=0$. Thus $M_{Y_2}\cong M_2$, leading to the contradiction as in the previous paragraph.

Finally, if\/ $B$ is of type (vi), then ${[M_{X_2}+M_{Y_2}:M_1]=0}$, $[M_{X_2}+M_{Y_2}:M_2]=1$ and ${[M_{X_2}+M_{Y_2}:M_2]=1}$, by \eqref{E:Smult} and Lemma \ref{L:MoritaTypes}. We may choose notation so that $M_{X_2}\cong M_2$ and $M_{Y_2}\cong M_3$. But $\ov{\chi_3}=\chi_4$. So $M_2^*\cong M_3$. We deduce that $M_{X_2}\cong M_{X_2}^*\cong M_{Y_2}$, which is a contradiction.
\end{proof}
\end{Theorem}

\section{Type {\rm (b)} extended defect groups}

In this section $E=D\langle e\rangle$, where $e^2=s_1$ and $e\in\op{C}(D)$. Table \ref{Tb:E} shows that the involutions in $E\backslash D$ form a single $E$-class, constisting of\/ $\{s_2e,s_2^{-1}e\}$. Theorem \ref{T:real_characters} shows that all height $1$ irreducible characters in $B$ are real. We require the following subsidiary result:

\begin{Lemma}\label{L:4Involutions}
Let $b$ be a real nilpotent $2$-block of a finite group $H$ that has defect couple $(Q,Q\times\langle r\rangle)$, where $Q$ is a Klein-four or dihedral group. Suppose that $\cclass_H(r)\cap Qr=\{r\}$. Then $k\cclass_H(r)b$ is irreducible.
\begin{proof}
We prove the result for the case that $Q\cong{\mathbb Z}_2^2$. The case that $Q$ is dihedral follows from similar arguments. Theorem 1.7 of \cite{Murray_GroupT} implies that $k\invol_Hb\cong S^4$, where $I$ is the unique irreducible $b$-module. So it is enough to show that $k\cclass_H(r)b$ is indecomposable. Then by the Green correspondence theorem, we may assume that $H=\op{N}_H(Q)$.

Let $b_1$ be a block of\/ $C:=\op{C}_H(Q)$ that is covered by $b$. Then $b_1$ is nilpotent, real, and has defect couple $(Q,Q\times\langle r\rangle)$. So $k\invol_Cb_1\cong S_1^4$, where $I_1$ is the unique irreducible $b_1$-module. The hypothesis implies that $Qr$ meets $4$ distinct conjugacy classes of involutions of\/ $C$. Using Theorem \ref{T:C_D(t)}, we deduce that $k\cclass_{C}(qr)b_1=I_1$, for each $q\in Q$.

As $\op{N}_H(b_1)=C$, we have $I=I_1{\uparrow^H}$ and hence $[I{\downarrow_{C}}:I_1]=1$. We then conclude from the previous paragraph that $k\cclass_{H}(r)b=I$.
\end{proof}
\end{Lemma}

Let $X/Y/Z$ denote a module with successive Loewy factors $X,Y,Z$.

\begin{Theorem}\label{T:vertexS}
$B$ is of type (i), (ii) or (v) and $k\invol B$ is indecomposable with vertex $S$. All $2^{d-3}$ irreducible characters in $F_{d-3}$ have FS-indicator $-1$ and the remaining $2^{d-3}+3$ irreducible characters in $B$ have FS-indicator $+1$. Moreover by the type of\/ $B$ we have:
\begin{itemize}
\item[(i)]  $k\invol B\cong M_1/M_1$.
\item[(ii)] $k\invol B\cong M_1/M_2/M_1$.
\item[(v)]  $k\invol B\cong M_1/(M_2\oplus M_3)/M_1$.
\end{itemize}
\begin{proof}
Suppose that $(t,b_t)\not\sim(s_1,b_{s_1})$. Then \cite{BrauerDihedral} shows that $b_t=b_{X_2}^{\op{C}(t)}$ and both $b_{X_2}$ and $b_t$ have defect group $X_2$. Lemma \ref{L:real_subpairs} implies that $b_{X_2}$ is real, with defect couple $(X_2,X_2\langle e\rangle)$. Then by the results in Section \ref{S:Introduction}, $b_t$ is real, and also has defect couple $(X_2,X_2\langle e\rangle)$. As $X_2\langle e\rangle$ does not split over $X_2$, it follows from Theorem \ref{T:C_D(t)} that $k\invol_{\op{C}(t)}b_t=0$. Similarly $(st,b_{st})\sim(s_1,b_{s_1})$ or $k\invol_{\op{C}(st)}b_{st}=0$.

Set ${H\!:=\!\op{C}(s_1)}$. Then the first paragraph implies that $b_{s_1}$ is the unique $2$-block of\/ $H$ such that $k\invol_Hb_{s_1}\ne0$ and $b_{s_1}^G=B$.

Table \ref{Tb:E} shows that there is a single $E$-class of involutions in $E\backslash D$, consisting of\/ $\{s_2e,s_2^{-1}e\}$. It then follows from Theorem \ref{T:C_D(t)} that $k\invol B=k\cclass(s_2e)B$. Moreover, as ${S=\op{C}_D(s_2e)}$, each component of\/ $k\invol B$ has vertex $S_i$, for some $i>0$, and at least one component has vertex $S$. As $\op{N}(S_i)\leq H$, Green correspondence induces a bijection between the components of\/ $k\invol B$ that have vertex $S_i$, and the components of\/ $k\cclass_H(s_2e)b_{s_1}$ that have vertex $S_i$.

Now $b_{s_1}$ is nilpotent and real, and has defect couple $(D,D\langle e\rangle)$. Let $I$ be the unique irreducible $b_{s_1}$-module and set $\ov{H}:=H/\langle s_1\rangle$.  Now $b_{s_1}$ dominates a unique block $\ov{b}$ of\/ $\ov{H}$, and $\ov{b}$ is nilpotent. Moreover, Lemma \ref{L:Dominates} implies that $\ov{b}$ is real, and has defect couple $(\ov{D},\ov{D}\times\langle\ov{e}\rangle)$. Table \ref{Tb:E} shows that $\ov{e},\ov{s_2e},\ov{te}$ and $\ov{ste}$ represent the distinct $\ov{E}$-classes of involutions in $\ov{E}\backslash\ov{D}$. Of these, only $\ov{s_2e}$ is the image of an involution in $H$. We conclude from Lemma \ref{L:4Involutions} that $k\cclass_{\ov{H}}(\ov{s_2e})\ov{b}\cong I$.

Now ${[\op{N}_H(\ov{s_2e})\!:\!\op{C}_H(s_2e)]=2}$, and $\op{N}_H(\ov{s_2e})$ is the preimage of\/ $\op{C}_{\ov{H}}(\ov{s_2e})$ in $H$. Then as $k\cclass_{\ov{H}}(\ov{s_2e})\ov{b}\cong I$, we get an exact sequence:
$$
0\rightarrow I\rightarrow k\invol_Hb\rightarrow I\rightarrow 0,\qquad\mbox{as $kH$-modules.}
$$
We cannot have $k\invol_Hb\cong I\oplus I$, as $I$ has no vertex contained in $S$. So $k\invol_Hb_{s_1}$ is indecomposable, with vertex $S$. We conclude from Lemma \ref{L:detectVertex} that $k\invol B$ is indecomposable, has vertex $S$, and is in Green correspondence with $k\invol_Hb_{s_1}$.

We have $[k\invol_Hb_{S_1}:I]=2$. Let $\theta$ be the Brauer character of\/ $I$. We temporarily apply the FS-indicator notation $\epsilon_i,\epsilon^{(j)}$ to $\op{Irr}(b_{s_1})$. Then from \eqref{E:Smult} and the decomposition matrix of\/ $b_{s_1}$ we get:
$$
\sum_{\psi\in\op{Irr}(b_{s_1})}\epsilon(\psi)d_{\psi,\theta}=\sum_{i=1}^4\epsilon_i+2\sum_{j=0}^{d-3}\epsilon^{(j)}2^j=2.
$$
The only solution is $\epsilon_i=1$, for $i=1,2,3,4$, $\epsilon^{(j)}=1$, for $j=0,\ldots,d-4$ and $\epsilon^{(d-3)}=-1$. In fact, all irreducible characters in $\ov{b}$ have FS-indicator $1$, and these account for all characters in $\op{Irr}(b_{s_1})$ that do not belong to the family $F_{d-3}$. We now use Lemma \ref{L:MoritaTypes} to compute
$$
\sum_{\psi\in\op{Irr}(b_{s_1})}\epsilon(\psi)d_{\psi,\theta}^{(s_1)}=\pm(4+2(2^{d-3}-1)-(-1)2.2^{d-3})=\pm(2^{d-1}+2).
$$
It then follows from Lemma \ref{L:2ElementSum}, and the fact that $s_1=e^2$, that
\begin{equation}\label{E:type(b)}
\sum_{\chi\in\op{Irr}(B)}\epsilon(\chi)d_{\chi,\theta}^{(s_1)}=2^{d-1}+2.
\end{equation}

We claim that $B$ is not of type (iii). For otherwise $\epsilon_i=+1$, for $i=1,2,3,4$, and we use the decomposition matrix of\/ $B$ to compute
$$
\sum_{\chi\in\op{Irr}(B)}\epsilon(\chi)d_{\chi\phi}^{(s_1)}=2\varepsilon\sum_{j=0}^{d-3}\epsilon^{(j)}2^j<2^{d-1}.
$$
This contradicts \eqref{E:type(b)}, proving our claim. Similar contradictions rule out the Morita types (iv) and (vi).

Now suppose that $B$ is of type (i), (ii) or (v). Then from the decomposition matrix of\/ $B$, and \eqref{E:type(b)} we get
$$
\varepsilon(\sum_{i=1}^4\epsilon_i+2\sum_{j=0}^{d-4}\epsilon^{(j)}2^j-2\epsilon^{(d-3)}2^{d-3})=\sum_{\chi\in\op{Irr}(B)}\epsilon(\chi)d_{\chi\phi}^{(s_1)}=2^{d-1}+2.
$$
In particular $\sum\epsilon_i\equiv0$ mod $4$. As $\epsilon_i\geq0$, we get $\epsilon_i=+1$, for $i=1,2,3,4$. Substituting these values into the equation above, we get
$$
\sum_{j=0}^{d-4}\epsilon^{(j)}2^j-\epsilon^{(d-3)2^{d-3}}=
\left\{
\begin{array}{rl}
 2^{d-2}-1,&\mbox{if\/ $\varepsilon=+1$;}\\
-2^{d-2}-3,&\mbox{if\/ $\varepsilon=-1$.}
\end{array}
\right.
$$
As the left hand side has absolute value $<2^{d-2}$, we deduce that $\varepsilon=+1$. The resulting equation has the unique solution $\epsilon^{(j)}=+1$, for $j=0,\ldots,d-4$ and $\epsilon^{(d-3)}=-1$.

Let $\theta_1$ be the Brauer character of\/ $M_1$. Then \eqref{E:Smult} gives 
$$
[k\invol:M_1]=\sum_{\chi\in\op{Irr}(B)}\epsilon(\chi)d_{\chi,\theta_1}=4+2(2^{d-3}-1)-2.2^{d-3}=2.
$$
It follows from this that $k\invol B=M_1/M_1$, if\/ $B$ is nilpotent.

Now suppose that $B$ has type (ii) or (v). Let $\theta_2$ be the Brauer character of\/ $M_2$. Using \eqref{E:Smult} and the decomposition matrix of\/ $B$, we have:
$$
[k\invol B:M_2]=\sum_{\chi\in\op{Irr}(B)}\epsilon(\chi)d_{\chi,\theta_2}=2+(2^{d-3}-1)-2^{d-3}=1.
$$
In the same way ${[k\invol B:M_3]\!=\!1}$, if\/ $B$ is of Morita type (v). We deduce the structure of the successive Loewy layers of\/ $k\invol B$ from these facts.
\end{proof}
\end{Theorem}

\section{Type {\rm(c)}: dihedral extended defect groups}

In this section $B$ is a real non-principal $2$-block with a dihedral defect group $D$ and a dihedral extended defect group $E=D\langle e\rangle$, where $e$ has order $2$ and $s=(te)^2$. In particular $B$ is not the principal $2$-block. We denote the projective cover of a module $M$ by $P(M)$.

By Theorem \ref{T:real_characters}, all height $1$ irreducible characters in $B$ are real-valued. Now $s_1$ is not the square of an element of\/ $E\backslash D$. It then follows from Lemma \ref{L:2ElementSum} that
\begin{equation}\label{E:s1_Dihedral}
\sum_{\chi\in\op{Irr}(B)}\epsilon(\chi)d_{\chi\theta}^{(s_1)}=0.
\end{equation}

\begin{Lemma}
$B$ is not of type (ii), (iii) or (iv).
\begin{proof}
Lemma \ref{L:cd13} implies that $B$ is not of type (ii) or (iii).

Suppose that $B$ is of type (iv). An examination of the decomposition matrix of\/ $B$ shows that all three irreducible $B$-modules are self-dual. It then follows from Theorem \ref{T:real_characters} that all irreducible characters in $B$ are real-valued. Corollary \ref{C:heightzero_FS+1} shows that $\epsilon_i=+1$, for $i=1,2,3,4$. We can now use the decompostion matrix of\/ $B$ to compute
$$
\sum_{\chi\in\op{Irr}(B)}\epsilon(\chi)d_{\chi\theta}^{(s_1)}=2\varepsilon(\sum_{j=0}^{d-4}\epsilon^{(j)}2^j-\epsilon^{(d-3)}2^{d-3})\equiv2\mod 4.
$$
This contradiction of \eqref{E:s1_Dihedral} completes the proof.
\end{proof}
\end{Lemma}

\begin{Theorem}\label{T:kOmega_dihedral}
$k\invol B$ is indecomposable, has vertex $\op{Z}(D)$, and is isomorphic to its Heller translate. Moreover, by the Morita type we have:
\begin{itemize}
\item[(i),(v)] Without loss of generality $\ov{\chi_3}=\chi_4$. All other irreducible\newline characters in $B$ have FS-indicator $+1$ and $P(k\invol B)\cong P(M_1)$;
\item[(vi)] All irreducible characters in $B$ have FS-indicator $+1$, and\newline $P(k\invol B)\cong P(M_1)\oplus P(M_2)\oplus P(M_3)$.
\end{itemize}
\begin{proof}
There is a single $E$-conjugacy class of involutions in $E\backslash D$, with representative $e$. Now $\op{C}_D(e)=\op{Z}(D)=\langle s_1\rangle$ is cyclic of order $2$. So by Theorem \ref{T:C_D(t)}, all components of\/ $k\invol B$ have vertex $\op{Z}(D)$.

Suppose first that $B$ is nilpotent. Lemmas \ref{L:real_subpairs} and \ref{L:subsection_conjugacy} show that $(t,b_t)$ and $(st,b_{st})$ are not real. Then by Theorem \ref{T:real_characters}, $B$ has two nonreal irreducible characters of height zero. We choose notation so that $\ov{\chi_3}=\chi_4$. Let $\theta_1$ be the Brauer character of\/ $M_1$ and let $M$ be a component of\/ $k\invol B$. Now $\op{dim}(M)=[M:M_1]\op{dim}(M_1)$, and $\nu\op{dim}(M_1)=\nu[G:D]$. So $\nu\op{dim}(M)=\nu[M:M_1]+\nu[G:D]$. As $M$ has vertex $\op{Z}(D)$, we have $\nu\op{dim}(M)\geq\nu[D:\op{Z}(D)]+\nu[G:D]$. We deduce that $2^{d-1}=[D:\op{Z}(D)]$ divides $[M:M_1]$. On the other hand, $[k\invol B:M_1]=\sum_{\chi\in\op{Irr}(B)}\nu(\chi)d_{\chi,\theta_1}$, which is $\leq2^{d-1}+2$, from the decomposition matrix of\/ $B$. We conclude that $[M:M_1]=2^{d-1}$, that $k\invol B=M$ is indecomposable, and that all irreducible characters in $B$, apart from $\chi_3,\chi_4$, have FS-indicator $+1$. Its easy to see that $k\invol B$ is the unique $B$-module that has vertex $\op{Z}(D)$. So $k\invol B$ coincides with its Heller translate $\Omega(k\invol B)$. Now from the Cartan matrix of\/ $B$, we have $[P(M_1):M_1]=2^d$. It follows that $P(k\invol B)\cong P(M_1)$.

Suppose that $B$ is of type (v) or (vi). Then $(t,b_t)\sim(s_1,b_{s_1})$ and $(st,b_{st})\sim(s_1,b_{s_1})$. So  $b_{s_1}$ is the unique $2$-block of\/ $H:=\op{C}(s_1)$ with $b_{s_1}^G=B$. The Green correspondence theorem establises a bijection between the components of\/ $k\invol B$ and the components of\/ $k\invol_Hb_{s_1}$. Applying the previous paragraph to $b_{s_1}$, we see that $k\invol_Hb_{s_1}$ is indecomposable and $k\invol_Hb_{s_1}=\Omega(k\invol_Hb_{s_1})$. So $k\invol B$ is indecomposable and $k\invol B=\Omega(k\invol B)$. We write $P(k\invol B)\cong P(M_1)^a\oplus P(M_2)^b\oplus P(M_3)^c$, where $a,b,c\geq0$.

We now specialize to the case that $B$ has Morita type (v). We see from the decomposition matrix of\/ $B$, and \eqref{E:s1_Dihedral} that
$$
\sum_{i=1}^4\epsilon_i+2\sum_{j=0}^{d-4}\epsilon^{(j)}2^j-2\epsilon^{(d-3)}2^{d-3}=0
$$
The only solution is $\epsilon_1=\epsilon_2=1$, $\ov{\chi_3}=\chi_4$, and $\epsilon^{(j)}=1$, for $j=0,\ldots,d-3$.  As $\Omega(k\invol B)=k\invol B$, we can compute $[k\invol B:M_i]$ in terms of the unknown exponents $a,b,c$, and compare these with the values given by \eqref{E:Smult}:
$$
\begin{array}{c|c|c|c}
  & M_1 & M_2 & M_3\\
\hline
\sum\epsilon(\chi)d_{\chi,M_i} & 2^{d-1} & 2^{d-2} & 2^{d-2}\\
P(\invol B) & 2^{d-1}(a+\!\frac{(b+c)}{2}) & 2^{d-2}(a+\!\frac{(b+c)}{2})\!+\!\frac{b}{2} & 2^{d-2}(a+\!\frac{(b+c)}{2})\!+\!\frac{c}{2}
\end{array}
$$
The only solution is $a=1,b=c=0$. So $P(k\invol B)\cong P(M_1)$.

Finally, we consider the case that $B$ is of Morita type (vi). We see from the decomposition matrix of\/ $B$, and \eqref{E:s1_Dihedral} that
$$
-\epsilon_1+\epsilon_2+\epsilon_3+\epsilon_4+2\sum_{j=0}^{d-4}\epsilon^{(j)}2^j-2\epsilon^{(d-3)}2^{d-3}=0
$$
In this case the only solution is $\epsilon_i=\epsilon^{(j)}=+1$, for all $i,j$. Proceeding as in the previous paragraph, we get the following table involving $a,b,c$:
$$
\begin{array}{c|c|c|c}
  & M_1 & M_2 & M_3\\
\hline
\sum\epsilon(\chi)d_{\chi,M_i} & 2 & 2^{d-2}+1 & 2^{d-2}+1\\
P(\invol B) & a+\frac{(b+c)}{2} & 2^{d-3}(b+c)+\frac{a+b}{2} & 2^{d-3}(b+c)+\frac{a+c}{2}
\end{array}
$$
Thus $a=b=c=1$, whence $P(k\invol B)\cong P(M_1)\oplus P(M_2)\oplus P(M_3)$.
\end{proof}
\end{Theorem}

\section{Type {\rm(d)}: semi-dihedral extended defect groups}

In this section $B$ is a real non-principal $2$-block with a dihedral defect group $D$ and a semi-dihedral extended defect group $E=D\langle e\rangle$, where $e^2=s_1$ and $s_1s=(te)^2$. By Theorem \ref{T:real_characters}, all height $1$ irreducible characters in $B$ are real-valued. As $E$ does not split over $D$, Theorem \ref{T:C_D(t)} implies that $k\invol B=0$. If\/ $M_i$ is an irreducible $B$-module, we use $\theta_i$ to denote its Brauer character. Then by \eqref{E:Smult} we have
\begin{equation}\label{E:zerosum}
\sum_{\chi\in\op{Irr}(B)}\epsilon(\chi)d_{\chi,\theta_i}=0.
\end{equation}

\begin{Theorem}
$B$ is of Morita type (i) or (v). We may choose notation so that $\ov{\chi_3}=\chi_4$; the irreducible characters in $F_{d-3}$ have FS-indicator $-1$; all other irreducible characters have FS-indicator $+1$.
\begin{proof}
Lemma \ref{L:cd13} implies that $B$ is not of type (ii) or (iii).

Suppose that $B$ is of Morita type (iv). Then from the decomposition matrix
$\sum_{\chi\in\op{Irr}(B)}\epsilon(\chi)d_{\chi,\theta_1}=\sum_{i=1}^4\epsilon_i>0$. This contradicts \eqref{E:zerosum}. So $B$ is not of type (iv). A similar argument shows that $B$ is not of type (vi).

Assume from now on that $B$ is of Morita type (i) or (v). Taking $i=1$ in \eqref{E:zerosum}, the decomposition matrix of\/ $B$ shows that
$$
\sum_{j=0}^{d-3}\epsilon^{(j)}2^j=-\sum_{i=1}^4\epsilon_i/2.
$$
Now the left hand side is odd and $\epsilon_i\geq 0$, for $i=1,2,3,4$. So without loss of generality $\epsilon_1=\epsilon_2=1$ and $\ov{\chi_3}=\chi_4$. We then get the unique solution $\epsilon^{(j)}=+1$, for $j=0,\ldots,d-4$ and $\epsilon^{(d-3)}=-1$.
\end{proof}
\end{Theorem}

\section{Type {\rm (e)} extended defect groups}\label{S:Type(e)}

In this section $B$ is a real non-principal $2$-block with extended defect group $E=D\langle e\rangle$, where $|D|\geq16$, $e^2=1$, $s^e=s_1s$ and $t^e=t$. By Theorem \ref{T:real_characters}, no character in $F_{d-3}$ is real, but the remaining characters of height $1$ in $F_0,\ldots,F_{d-4}$ are real.

\begin{Theorem}\label{T:type(e)}
$B$ is of type (i), (ii) or (v). Let $\chi\in\op{Irr}(B)$. Then $\chi$ is nonreal, if\/ $\chi\in F_{d-3}$. Otherwise $\epsilon(\chi)=+1$. We may write ${k\invol B=M_{X_{d-1}}\oplus M_{X_2}}$, where $M_{X_i}$ is an indecomposable module with vertex $X_i$, for $i=d-1,2$. If\/ $B$ is nilpotent then $[M_{X_i}:M_1]=2^{d-i}$.
\begin{proof}
By Table \ref{Tb:E}, there are two $E$-conjugacy classes of involutions in $E\backslash D$, with representatives $e$ and $te$. Now $\op{C}_D(e)=X_{d-1}$ and $\op{C}_D(te)=X_2$. So part (ii) of Theorem \ref{T:C_D(t)} implies that $k\cclass_G(e)B$ has at least one component $M_{X_{d-1}}$ that has vertex $X_{d-1}$.

Set $b_2:=b_{X_2}^{\op{N}(X_2)}$. Then $b_2$ is real and nilpotent, and has defect couple $(X_3,X_3\times\langle e\rangle)$. As $X_3\cong D_8$, Theorem \ref{T:TotallySplitModules} implies that $k\invol_{\op{N}(X_2)}b_2$ has one component with vertex $X_2$. Set $M_{X_2}$ as its Green correspondent with respect to $(G,X_2,\op{N}(X_2))$. Then $M_{X_2}$ is the unique component of\/ $k\invol B$ that has vertex $B$-subpair $(X_2,b_{X_2})$.

Suppose that $B$ is nilpotent. This holds for example if\/ $G=\op{C}(s_1)$. Then $B$ is real and has defect couple $(D,E)$. Let $\Phi$ be the unique principal indecomposable character of\/ $B$. Then from \eqref{E:Smult} and the decomposition matrix of\/ $B$, we have $[k\invol:M_1]=\epsilon(\Phi)\leq2^{d-2}+2$. On the other hand, $2=[D:X_{d-1}]$ divides $[M_{X_{d-1}}:M_1]$, and $2^{d-2}=[D:X_2]$ divides $[M_{X_2}:M_1]$. We deduce that $k\invol B=M_{X_{d-1}}\oplus M_{X_2}$. Moreover, $\epsilon(\Phi)=2^{d-2}+2$. So $\epsilon(\chi)=+1$, if\/ $\chi$ is an irreducible character in $B$ that does not belong to $F_{d-3}$. The last statement of the Theorem also follows.

We now remove the nilpotency assumption on $B$. The last paragraph applies to $b_{s_1}$, allowing us to compute
\begin{equation}\label{E:type(e)}
\sum_{\psi\in\op{Irr}(b_{s_1})}\epsilon(\psi)d_{\psi,\theta}^{(s_1)}=2^{d-2}+2.
\end{equation}
Set $b_{d-1}:=b_{X_{d-1}}^{\op{N}(X_{d-1})}$. Then $b_{d-1}$ is real and nilpotent, and has defect couple $(D,E)$. The last paragraph shows that $k\invol_{\op{N}(X_{d-1})}b_{d-1}$ has one component with vertex $X_{d-1}$. Set $M_{X_{d-1}}$ as its Green correspondent with respect to $(G,X_{d-1},\op{N}(X_{d-1}))$. Then $M_{X_{d-1}}$ is the unique component of\/ $k\invol B$ that has vertex $B$-subpair $(X_{d-1},b_{X_{d-1}})$.

Now let $(V,b_V)$ be a vertex $B$-subpair of a component of\/ $k\invol B$. Lemmas \ref{L:detectVertex} and \ref{L:subpair_conjugacy} imply that $(V,b_V)$ is a vertex $B$-subpair of a component of\/ $k\invol_{\op{C}(x)}b_x$, where $x=s_1,t$ or $st$. Now $b_{st}$ is real with vertex pair $(Y_2,Y_2\langle s_2e\rangle)$. As $Y_2\langle s_2e\rangle$ does not split over $Y_2$, Theorem \ref{T:C_D(t)} implies that $k\invol_{\op{C}(st)}b_{st}=0$. Our work above shows that every component of\/ $k\invol_{\op{C}(s_1)}b_{s_1}$ has vertex $B$-subpair $(X_i,b_{X_i})$, for $i={d-1},2$. Theorem 1.7(i) of \cite{Murray_GroupT} shows that every component of\/ $k\invol_{\op{C}(t)}b_t$ has vertex $B$-subpair $(X_2,b_{X_2})$. We deduce from this that $k\invol B=M_{X_{d-1}}\oplus M_{X_2}$. This proves one statement of the Theorem.

Suppose for the sake of contradiction that $B$ is of type (iii). Then all irreducible characters of height $0$ are real, whence they have FS-indicator $+1$. From the decomposition matrix of\/ $B$, we compute
$$
\sum_{\chi\in\op{Irr}(B)}\epsilon(\chi)d_{\chi\phi}^{(s_1)}=2\sum_{j=0}^{d-4}\epsilon^{(j)}2^j<2^{d-2}.
$$
This contradicts Lemma \ref{L:BrauerSum} and \eqref{E:type(e)}. We can show that $B$ is not of Morita type (iv) or (vi), using similar arguments.

Now suppose that $B$ is of type (ii) or (v). Then we compute
$$
\varepsilon(\sum_{i=1}^4\epsilon_i+2\sum_{j=0}^{d-4}\epsilon^{(j)}2^j)=\sum_{\chi\in\op{Irr}(B)}\epsilon(\chi)d_{\chi\theta}^{(s_1)}=2^{d-2}+2,
$$
using Lemma \ref{L:BrauerSum} and \eqref{E:type(e)}. Considering this equality modulo $4$, and using the fact that $\epsilon_i\in\{0,1\}$, we see that $\epsilon_i=+1$, for $i=1,2,3,4$. Substituting these values, we get two possibilities
$$
\sum_{j=0}^{d-4}\epsilon^{(j)}2^j=
\left\{
\begin{array}{rl}
 2^{d-3}-1,&\mbox{if\/ $\varepsilon=+1$;}\\
-2^{d-3}-3,&\mbox{if\/ $\varepsilon=-1$.}
\end{array}
\right.
$$
As the left hand side has absolute value $<2^{d-3}$, we deduce that $\varepsilon=+1$. The resulting equation has the unique solution $\epsilon^{(j)}=+1$, for $j=0,\ldots,d-4$. This completes the proof.
\end{proof}
\end{Theorem}

\section{Summary of results for blocks with dihedral defect groups}

We summarize the results of Sections \ref{S:TotallySplit} through \ref{S:Type(e)} in the following table. Recall that $M_V$ is an indecomposable module that has vertex $V$.
\begin{table}[h]\begin{center}
\begin{tabular}{|c|c|c|c|c|c|c}
\hline
Block& $E$  & $\epsilon_i$    & $\epsilon^{(j)}$  & $\!\epsilon^{(d-3)}\!$ & $k\invol B$\\
Type & type & ${\scriptstyle 1\leq i\leq 4}$ & ${\scriptstyle 0\leq j\leq d-4}$ &  & \\ \hline\hline
 & (a) & $++++$ & $+$ & $+$ & $M_1^2\oplus M_{X_2}\oplus M_{Y_2}$\\
          & (b) & $++++$ & $+$ & $-$ & $M_S=M_1/M_1$\\
 Nilpotent         & (c) & $++00$ & $+$ & $+$ & $M_{\mathord{Z}(D)}$\\
          & (d) & $++00$ & $+$ & $-$ &  -- \\
          & (e) & $++++$ & $+$ & $0$ & $M_{X_{d-1}}\oplus M_{X_2}$\\ \hline\hline
$\op{PGL}(2,q)$ & (a) & $++++$ & $+$ & $+$ & $M_D^2\oplus M_{X_2}\oplus M_{Y_2}$\\
$q\equiv1\mbox{(mod $4$)}$& (b) & $++++$ & $+$ & $-$ & $M_S=M_1/M_2/M_1$\\
          & (e) & $++++$ & $+$ & $0$ & $M_{X_{d-1}}\oplus M_{X_2}$\\ \hline\hline
$\op{PGL}(2,q)$ & (a) & $++++$ & $+$ & $+$ & $M_D^2\oplus M_{X_2}\oplus M_{Y_2}$\\
$q\equiv3\mbox{(mod $4$)}$&     &        &     &     & \\ \hline\hline
$A_7$     & (a) & $++++$ & $+$ & $+$ & $M_D^2\oplus M_{X_2}\oplus M_{Y_2}$\\ \hline\hline
 & (a) & $++++$ & $+$ & $+$ & $M_D^2\oplus M_{X_2}\oplus M_{Y_2}$\\
$\op{PSL}(2,q)$& (b) & $++++$ & $+$ & $-$ & $\!M_S\!=\!M_1\!/\!(\!M_2\!\oplus\! M_3\!)\!/\!M_1\!\!$\\
$q\equiv1\mbox{(mod $4$)}$          & (c) & $++00$ & $+$ & $+$ & $M_{\mathord{Z}(D)}$\\
          & (d) & $++00$ & $+$ & $-$ &  -- \\
          & (e) & $++++$ & $+$ & $0$ & $M_{X_{d-1}}\oplus M_{X_2}$\\ \hline\hline
$\op{PSL}(2,q)$ & (a) & $++00$ & $+$ & $+$ & $M_D^2\oplus M_{X_2}\oplus M_{Y_2}$\\
$q\equiv3\mbox{(mod $4$)}$& (c) & $++++$ & $+$ & $+$ & $M_{\mathord{Z}(D)}$\\ \hline
\end{tabular}
\end{center}\caption{Summary of results}\label{Tb:main}\end{table}

\newpage

\end{document}